\documentclass[12pt]{article}
\usepackage{amssymb}
\newcommand{\R}{\mathbb{R}}

\newcommand{\Z}{\mathbb{Z}}
\renewcommand{\P}{\mathbb{P}}
\newcommand{\Q}{\mathbb{Q}}
\newcommand{\E}{\mathbb{E}}

\newcommand{\N}{\mathbb{N}}

\newcommand{\T}{\mathbb{T}}
\def\build#1_#2^#3{\mathrel{
\mathop{\kern 0pt#1}\limits_{#2}^{#3}}}
\def\llbracket{[\hspace{-.10em} [ }
\def\rrbracket{ ] \hspace{-.10em}]}
\def\cq{$\hfill \square$}
\def \un{\underline}
\def\ind{{\bf 1}_}

\def\vep{\varepsilon}

\def\sg{\sigma}
\def\dl{\delta}
\def\l{{\cal L}}

\def\m{{\cal M}}
\def\n{{\cal N}}

\def\e{{\cal E}}

\def\u{{\cal U}}

\def\t{{\cal T}}

\def\r{{\cal R}}
\def\p{{\cal P}}
\def\q{{\cal Q}}

\def\z{{\cal Z}}

\def\W{{\cal W}}
\def\w{{\rm w}}

\def\be{\begin{equation}}
\def\ee{\end{equation}}
\def\ba{\begin{eqnarray*}}
\def\ea{\end{eqnarray*}}
\def\ov{\overline}
\def\wh{\widehat}
\def\wt{\widetilde}
\def\cg{\Big[}
\def\cd{\Big]}
\def\pg{\Big(}
\def\pd{\Big)}
\def\la{\longrightarrow}
\def\da{\downarrow}

\def\proof{\vskip 3mm \noindent{\bf Proof:}\hskip10pt}

\newtheorem{theorem}{Theorem}[section]
\newtheorem{lemma}[theorem]{Lemma}
\newtheorem{proposition}[theorem]{Proposition}
\newtheorem{corollary}[theorem]{Corollary}

\begin{document}

\title{ \bf CONDITIONED BROWNIAN TREES}
\author{
Jean-Fran\c cois {\sc Le Gall} and Mathilde {\sc Weill} \\
{\small D.M.A., Ecole normale sup\'erieure, 45 rue d'Ulm, 75005 Paris, France}}
\vspace{2mm}
\date{\today} 

\maketitle

\begin{abstract}
We consider a Brownian tree consisting of a collection of one-dimensional
Brownian paths started from the origin, whose genealogical structure is given by
the Continuum Random Tree (CRT). This Brownian tree may be generated 
from the Brownian snake driven by a normalized Brownian 
excursion, and thus yields a convenient representation 
of the so-called Integrated Super-Brownian Excursion (ISE), 
which can be viewed as the uniform probability measure on the tree of paths. We discuss
different approaches that lead to the definition of the Brownian tree conditioned to stay on
the positive half-line. We also establish a Verwaat-like theorem showing that this
conditioned Brownian tree can be obtained by re-rooting the unconditioned one
at the vertex corresponding to the minimal spatial position. In terms of ISE, this theorem yields
the following fact: Conditioning ISE to put no mass 
on $]-\infty,-\varepsilon[$ and letting $\varepsilon$ go to $0$ is equivalent to
shifting the unconditioned ISE to the right so that the left-most point of its
support becomes the origin.  We derive a number of explicit estimates and formulas
for our conditioned Brownian trees. In particular, the probability that ISE puts no
mass on $]-\infty,-\varepsilon[$ is shown to behave like $2\varepsilon^4/21$
when $\varepsilon$ goes to $0$. Finally, for the conditioned Brownian tree with a fixed
height $h$, we obtain a decomposition involving a spine whose distribution is absolutely
continuous with respect to that of a nine-dimensional Bessel process
on the time interval $[0,h]$, and Poisson
processes of subtrees originating from this spine.
\end{abstract}

\section{Introduction}

In this work, we define and study a continuous tree of one-dimensional Brownian paths
started from the origin, which is conditioned to remain in the positive half-line. An
important motivation for introducing this object comes from its relation with
analogous discrete models which are discussed in several recent papers.

In order to present our main results, let us briefly describe a construction
of unconditioned Brownian trees. We start from a positive Brownian 
excursion conditioned to have duration $1$ (a normalized Brownian excursion
in short), which is denoted by $(e(s),0\leq s\leq 1)$. This random function
can be viewed as coding a continuous tree via the following simple prescriptions.
For every $s,s'\in[0,1]$, we set
$$m_e(s,s'):=\inf_{s\wedge s'\leq r\leq s\vee s'}e(r).$$
We then define an equivalence relation on $[0,1]$ by setting $s\sim s'$ if and 
only if $e(s)=e(s')=m_e(s,s')$. Finally we put
$$d_e(s,s')=e(s)+e(s')-2\,m_e(s,s')$$
and note that $d_e(s,s')$ only depends on the equivalence classes of $s$
and $s'$. Then the quotient space ${\bf T}_e:=[0,1]/\sim$ equipped with the metric
$d_e$ is a compact $\R$-tree (see e.g. Section 2 of \cite{DuLG}). In other words, it is a compact metric space
such that for any two points $\sigma$ and $\sigma'$
there is a unique arc with endpoints $\sigma$ and $\sigma'$ and furthermore
this arc is isometric to a compact interval of the real line.
We view ${\bf T}_e$ as a rooted $\R$-tree, whose root $\rho$ is the equivalence
class of $0$. For every $\sigma\in{\bf T}_e$, the ancestral line 
of $\sigma$ is the line segment joining $\rho$ to $\sigma$. This line segment is denoted
by $\llbracket \rho,\sigma\rrbracket$. We write $\dot s$
for the equivalence class of $s$, which is a vertex in ${\bf T}_e$
at generation $e(s)=d_e(0,s)$.

Up to unimportant scaling constants, ${\bf T}_e$ is the Continuum Random Tree
(CRT) introduced by Aldous \cite{Al1}. The preceding presentation
is indeed a reformulation of Corollary 22 in \cite{Al3}, which was proved 
via a discrete approximation (a more direct approach was given in \cite{LG1}).
As Aldous \cite{Al3} has shown, the CRT is the scaling limit of critical Galton-Watson trees
conditioned to have a large fixed progeny (see \cite{Duq} and
\cite{DuLG} for recent generalizations of Aldous' result).
The fact that Brownian excursions can be used to model continuous genealogies had
been used before, in particular in the Brownian snake approach to
superprocesses (see \cite{LG0}).

We can now combine the branching structure of the CRT with independent
spatial motions. We restrict ourselves to spatial displacements given by
linear Brownian motions, which is the case of interest in this work. 
Conditionally given $e$, we introduce a centered Gaussian process
$(V_\sigma,\sigma\in{\bf T}_e)$ with covariance
$${\rm cov}(V_{\dot s},V_{\dot s'})=m_e(s,s')\ ,\qquad s,s'\in[0,1].$$
This definition should become clear if we observe that $m_e(s,s')$ is the
generation of the most recent common ancestor to $\dot s$ and $\dot s'$ in the tree
${\bf T}_e$. It is easy to verify that the process
$(V_\sigma,\sigma\in{\bf T}_e)$ has a continuous modification.
The random measure $\z$ on $\R$ defined by
$$\langle \z,\varphi\rangle =\int_0^1 \varphi(V_{\dot s})\,ds$$
is then the one-dimensional Integrated Super-Brownian Excursion (ISE). 
Note that ISE in higher dimensions, and related Brownian trees, have appeared recently in 
various asymptotic results for statistical mechanics models
(see e.g. \cite{Sl1},\cite{Sl2},\cite{Sl3}). The support, or range, of ISE is
$$\r:=\{V_\sigma:\sigma\in{\bf T}_e\}.$$

For our purposes, it is also convenient to reinterpret the preceding notions
in terms of the Brownian snake. The Brownian snake $(W_s,0\leq s\leq 1)$ driven
by the normalized excursion $e$ is obtained as follows
(see subsection 2.1 for
a more detailed presentation). For every 
$s\in[0,1]$, $W_s=(W_s(t),0\leq t\leq e(s))$ is the finite path which gives
the spatial positions along the ancestral line of $\dot s$: $W_s(t)=V_\sigma$
if $\sigma$ is the vertex at distance $t$ from the root on the segment
$\llbracket \rho,\dot s\rrbracket$. Note that $W_s$
only depends on the equivalent class $\dot s$. We view $W_s$ as a random element 
of the space $\W$ of finite paths.

Our first goal is to give a precise definition of the Brownian tree
$(V_\sigma,\sigma\in{\bf T}_e)$ conditioned to remain positive. Equivalently 
this amounts to conditioning ISE to put no mass on the negative half-line. Our first
theorem gives a precise meaning to this conditioning in terms of the
Brownian snake. We denote by $\N^{(1)}_0$ the distribution of
$(W_s,0\leq s\leq 1)$ on the canonical space $C([0,1],\W)$ of continuous
functions from $[0,1]$ into $\W$, and we abuse notation by still writing 
$(W_s,0\leq s\leq 1)$ for the canonical process on this space. The
range $\r$ is then defined under $\N^{(1)}_0$ by
$$\r=\{\wh W_s:0\leq s\leq 1\}$$
where $\wh W_s$ denotes the endpoint of the path $W_s$.

\begin{theorem}
\label{main-cond}
We have 
$$\lim_{\varepsilon\da
0}\varepsilon^{-4}\,\N^{(1)}_0(\r\subset]-\varepsilon,\infty[)={2\over 21}.$$
There
exists a probability measure on $C([0,1],\W)$, which is denoted by $\ov\N^{(1)}_0$,
such that
$$\lim_{\varepsilon\da 0} \N^{(1)}_0(\cdot\mid
\r\subset]-\varepsilon,\infty[)=\ov\N^{(1)}_0,$$ in the sense of weak convergence in
the space of probability measures on
$C([0,1],\W)$.
\end{theorem}

Our second theorem gives an explicit representation of the
conditioned measures $\ov\N^{(1)}_0$, which is analogous to a 
famous theorem of Verwaat \cite{Verwaat} relating the normalized Brownian excursion
to the Brownian bridge. To state this result, we need the notion of re-rooting.
For $s\in[0,1]$, we write ${\bf T}_e^{[s]}$ for the ``same'' tree ${\bf T}_e$
but with root $\dot s$ instead of $\rho=\dot 0$. We then shift the spatial 
positions by setting $V^{[s]}_\sigma=V_\sigma-V_{\dot s}$ for every $\sigma
\in{\bf T}_e$, in such a way that the spatial position of the new root is still the origin. 
(Notice that both  ${\bf T}_e^{[s]}$ and $V^{[s]}$ only depend on
$\dot s$, and we could as well define ${\bf T}_e^{[\sigma]}$ and $V^{[\sigma]}$
for $\sigma\in{\bf T}_e$.)
Finally, the re-rooted snake $W^{[s]}=(W^{[s]}_r,0\leq r\leq 1)$ is defined
analogously as before: For every $r\in[0,1]$, $W^{[s]}_r$ is the path
giving the spatial positions $V^{[s]}_\sigma$ along the ancestral line (in the
re-rooted tree) of the vertex $s+r$ mod. $1$.

\begin{theorem}
\label{verwaat}
Let $s_*$ be the unique time of the minimum of $\wh W$ on $[0,1]$.
The probability measure $\ov\N^{(1)}_0$ is the law under 
$\N^{(1)}_0$ of the re-rooted snake $W^{[s_*]}$.
\end{theorem}

If we want to define
one-dimensional ISE conditioned to put no mass on the negative half-line,
the most natural way is to condition it to put no mass on $]-\infty,-\varepsilon[$ 
and then to let $\varepsilon$ go to $0$. As a consequence of 
the previous two theorems, this is equivalent to shifting
the unconditioned ISE to the right, so that the left-most point
of its support becomes the origin. 

Both Theorem \ref{main-cond} and Theorem \ref{verwaat} could be presented 
in a different and perhaps more elegant manner by using the formalism of
spatial trees as in Section 5 of \cite{DuLG}. In this formalism, a spatial
tree is a pair $({\bf T},U)$ where ${\bf T}$ is a compact rooted $\R$-tree
(in fact an equivalent class of such objects modulo root-preserving
isometries) and $U$ is a continuous mapping from
${\bf T}$ into $\R^d$. Then the second 
assertion of Theorem \ref{main-cond} can be rephrased by saying that
the conditional distribution of the spatial tree $({\bf T}_e,V)$
knowing that $\r\subset]-\varepsilon,\infty[$ has a limit when $\varepsilon$ goes to $0$,
and Theorem \ref{verwaat} says that this limit is the distribution
of $({\bf T}^{[\sigma_*]}_e,V^{[\sigma_*]})$ where $\sigma_*$ is the unique vertex
minimizing $V$. We have chosen the above presentation because the Brownian snake
plays a fundamental role in our proofs and also because the resulting statements
are stronger than the ones in terms of spatial trees.

Let us discuss the relationship of the above theorems with previous results. The first assertion
of Theorem \ref{main-cond} is closely related to some
estimates of Abraham and Werner \cite{AW}. In particular, Abraham and Werner proved that the
probability for a Brownian snake driven by a Brownian excursion
of height $1$ not to hit the set $]-\infty,-\varepsilon[$ behaves like
a constant times $\varepsilon^4$ (see Section 4 below).
The $d$-dimensional  Brownian snake conditioned not to exit a domain
$D$ was studied by Abraham and Serlet \cite{AS}, who observed that this conditioning gives
rise to a particular instance of the Brownian snake with drift. The setting
in \cite{AS} is different from the present work, in that the initial point
of the snake lies inside the domain, and not at its boundary as here. We also mention
the paper \cite{JR} by Jansons and Rogers, who establish a decomposition at the
minimum for a Brownian tree where branchings occur only at discrete times.

An important motivation for the present work came from several recent
papers that discuss asymptotics for planar maps. A key result due to
Schaeffer (see \cite{CS}) establishes a bijection between rooted planar quadrangulations
and certain discrete trees called well-labelled trees. Roughly, a 
well-labelled tree consists of a (discrete) plane tree whose vertices are 
given labels which are positive integers, with the constraints that
the label of the root is $1$ and 
the labels of two neighboring vertices can differ by at most $1$. 
Our conditioned Brownian snake should then be viewed as a continuous
model for well-labelled trees. This idea was exploited in  \cite{CS}
and especially in Marckert and Mokkadem \cite{MM}, where the re-rooted snake $W^{[s_*]}$ appears 
in the description of the Brownian map, which is the continuous
object describing scaling limits of planar quadrangulations. In contrast 
with the present work, the re-rooted snake $W^{[s_*]}$ is not interpreted in \cite{MM} as
a conditioned object, but rather as a scaling limit of
re-rooted discrete snakes. Closely related models of discrete labelled trees are
also of interest in theoretical physics: See in particular \cite{BDG1} and \cite{BDG2}. Motivated by \cite{CS}
and \cite{MM}, we prove in \cite{LG} that our conditioned Brownian tree is
the scaling limit of discrete spatial trees conditioned to remain positive. To be specific, we consider
a  Galton-Watson tree whose offspring
distribution is critical and has (small) exponential moments, and we condition this tree
to have exactly $n$ vertices (in the special case of the 
geometric distribution, this gives rise to a tree that is uniformly
distributed over the set of plane trees with $n$ vertices).
This branching structure is combined with a spatial displacement 
which is a symmetric random walk with bounded jump size on $\Z$. 
Assuming that the root is at the origin of $\Z$, 
the spatial tree is then conditioned to remain on the positive side.
According to the main theorem of \cite{LG}, the scaling limit of this 
conditioned discrete tree when $n\to \infty$ leads to
the measure $\ov\N^{(1)}_0$ discussed above. The convergence here, and the 
precise form of the scaling transformation, are as in 
Theorem 2 of \cite{JM}, which discusses
scaling limits for unconditioned discrete snakes.

Let us now describe the other contributions of this paper.
Although the preceding theorems have been stated for the measure $\N^{(1)}_0$, 
a more fundamental object is the excursion measure $\N_0$ of the Brownian snake
(see e.g. \cite{Zu}). Roughly speaking, $\N_0$ is obtained by the same
construction as above, but instead of considering a normalized Brownian excursion, 
we now let $e$ be distributed according to the (infinite) It\^o measure 
of Brownian excursions. If $\sigma(e)$ denotes the duration of excursion $e$,
we have $\N^{(1)}_0=\N_0(\cdot\mid \sigma=1)$. It turns out that 
many calculations are more tractable under the infinite measure 
$\N_0$ than under $\N^{(1)}_0$.  For this reason, both Theorems \ref{main-cond} and Theorem
\ref{verwaat}  are proved in Section 3 as consequences of Theorem \ref{exc-cond}, which deals
with $\N_0$. Motivated by Theorem \ref{exc-cond} we introduce another infinite measure
denoted by $\ov\N_0$, which should be interpreted as $\N_0$ conditioned
on the event $\{\r\subset [0,\infty[\}$, even though the conditioning
requires some care as we are dealing with infinite measures. In the same
way as for unconditioned measures, we have $\ov\N^{(1)}_0=\ov\N_0(\cdot\mid \sigma=1)$.
Another motivation for considering the measure $\ov\N_0$ comes from
connections with superprocesses: Analogously to Chapter IV of \cite{Zu}
in the unconditioned case, $\ov\N_0$
could be used to define and to analyse a 
one-dimensional super-Brownian motion started from the Dirac measure $\delta_0$ and
conditioned never to charge the negative half-line.

In Section 4, we present a different approach that leads to the same limiting measures. 
If $H(e)$ stands for
the height of excursion $e$, we consider
for every $h>0$ the measure $\N^h_0:=\N_0(\cdot\mid H=h)$. In the above construction
this amounts to replacing the normalized excursion $e$ by a Brownian excursion with height $h$. 
By using a famous decomposition theorem of Williams, we can then analyse the 
behavior of the measure $\N^h_0$ conditioned on the event that the
range does not intersect $]-\infty,-\varepsilon[$ and show that it has 
a limit denoted by $\ov\N^h_0$ when $\varepsilon\to 0$. The method also provides
information about the Brownian tree under $\ov\N^h_0$: This Brownian tree
consists of a spine whose distribution is absolutely continuous with respect
to that of the nine-dimensional Bessel process, and as usual a Poisson collection
of subtrees originating from the spine, which are Brownian snake excursions 
conditioned not to hit the negative half-line. The connection
with the measures $\ov\N^{(1)}_0$ and $\ov\N_0$
is made by proving that $\ov\N^h_0=\ov\N_0(\cdot\mid H=h)$. Several arguments in this section have been 
inspired by Abraham and Werner's paper \cite{AW}. It should also be noted that a discrete 
version of the nine-dimensional Bessel process already appears in the 
Chassaing-Durhuus paper \cite{CD}.

At the end of Section 4, we also discuss the limiting behavior of the measures $\ov\N^h_0$
as $h\to \infty$. This leads to a probability measure $\ov\N^\infty_0$
that should be viewed as the law of an infinite Brownian snake excursion conditioned
to stay positive. We again get a description of the Brownian tree coded by
$\ov\N^\infty_0$ in terms of a spine and conditioned Brownian snake excursions
originating from this spine. Moreover, the description is simpler in the sense that
the spine is exactly distributed as a nine-dimensional Bessel process started at the
origin.

Section 5 gives an explicit formula for the finite-dimensional marginal 
distributions of the Brownian tree under $\ov\N_0$, that is for
$$\ov\N_0\Big(\int_{]0,\sigma[^p} ds_1\ldots
ds_p\,F(W_{s_1},\ldots,W_{s_p})\Big) $$
where $p\geq 1$ is an integer and $F$ is a symmetric nonnegative 
measurable function on $\W^p$. In a way similar to the corresponding
result for the unconditioned Brownian snake (see (\ref{margi}) below),
this formula involves combining the branching structure of certain discrete trees
with spatial displacements. Here however because of the 
conditioning, the spatial displacements turn out
to be given by nine-dimensional Bessel processes rather than linear Brownian 
motions. In the same way as the finite-dimensional marginal distributions of the CRT can be
derived from the analogous formula under the It\^o measure
(see Chapter III of \cite{Zu}), one might hope to derive the
expression of the finite-dimensional marginals under $\ov\N^{(1)}_0$
from the case of $\ov\N_0$. This idea apparently leads to untractable calculations, 
but we still expect Theorem \ref{marginalscond} to have useful applications
in future work about conditioned trees.

Basic facts about the Brownian snake are recalled in Section 2, 
which also establishes a few important preliminary results, some
of which are of independent interest. In particular,
we state and prove a general version of the invariance 
property of $\N_0$ under re-rooting (Theorem \ref{reroot}).
This result is clearly related to the invariance of the CRT
under uniform re-rooting, which was observed by Aldous \cite{Al2}
(and generalized to L\'evy trees in Proposition 4.8 of \cite{DuLG}).
See also \cite{CS} for similar ideas in a discrete setting, and
especially Proposition 13 of \cite{MM} which gives a closely
related statement.

\section{Preliminaries}

In this section, we recall the basic facts about the Brownian snake
that we will use later, and we also establish a few important
preliminary results. We refer to \cite{Zu} for a more detailed 
presentation of the Brownian snake and its connections with partial
differential equations. In the first four subsections below, we deal with
the $d$-dimensional Brownian snake since the proofs are not more
difficult in that case, and the results may have other applications.

\subsection{The Brownian snake}

The ($d$-dimensional) Brownian snake is a Markov process taking values
in the space $\W$ of finite paths in $\R^d$. Here a finite path is simply 
a continuous mapping $\w:[0,\zeta]\la \R^d$, where
$\zeta=\zeta_{(\w)}$ is a nonnegative real number called the 
lifetime of $\w$. The set $\W$ is a Polish space when equipped with the
distance
$$d(\w,\w')=|\zeta_{(\w)}-\zeta_{(\w')}|+\sup_{t\geq 0}|\w(t\wedge
\zeta_{(\w)})-\w'(t\wedge\zeta_{(\w')})|.$$
The endpoint (or tip) of the path $\w$ is denoted by $\wh \w$.
The range of $\w$ is denoted by $\w[0,\zeta_{(\w)}]$.

In this work, it will be convenient to use the
canonical space $\Omega:=C(\R_+,\W)$ of continuous functions from
$\R_+$ into $\W$, which is equipped with the topology
of uniform convergence on every compact subset of $\R_+$. The canonical process on $\Omega$ is then denoted by
$$W_s(\omega)=\omega(s)\;,\quad \omega\in\Omega\;,$$
and we write $\zeta_s=\zeta_{(W_s)}$ for the lifetime of $W_s$.

Let $\w\in\W$. The law of the Brownian snake started from $\w$
is the probability measure $\P_\w$ on $\Omega$ which can be 
characterized as follows. First, the process $(\zeta_s)_{s\geq 0}$
is under $\P_\w$ a reflected  Brownian motion in $[0,\infty[$ started 
from $\zeta_{(\w)}$. Secondly, the conditional distribution 
of $(W_s)_{s\geq 0}$ knowing $(\zeta_{s})_{s\geq 0}$, which is
denoted by $\Theta^\zeta_\w$, is characterized by the
following properties:
\begin{description}
\item{(i)} $W_0=\w$, $\Theta^\zeta_\w$ a.s.
\item{(ii)} The process $(W_s)_{s\geq 0}$ is time-inhomogeneous
Markov under $\Theta^\zeta_\w$. Moreover, if $0\leq s\leq
s'$,
\begin{description}
\item{$\bullet$} $W_{s'}(t)=W_{s}(t)$ for every $t\leq m(s,s'):=
\inf_{[s,s']}\zeta_r$, \ $\Theta^\zeta_\w$ a.s.
\item{$\bullet$} $(W_{s'}(m(s,s')+t)-W_{s'}(m(s,s')))_{0\leq t\leq
\zeta_{s'}- m(s,s')}$ is independent of $W_s$ and distributed as a
$d$-dimensional Brownian motion started at $0$  under
$\Theta^\zeta_\w$.
\end{description}
\end{description}
Informally, the value $W_s$ of the Brownian snake at time $s$
is a random path with a random lifetime $\zeta_s$ evolving like
reflecting Brownian motion in $[0,\infty[$. When $\zeta_s$ decreases,
the path is erased from its tip, and when $\zeta_s$ increases, the path 
is extended by adding ``little pieces'' of Brownian paths at its tip.

Excursion measures play a fundamental role throughout this work. 
We denote by $n(de)$ the It\^o measure of positive Brownian excursions.
This is a $\sigma$-finite measure on the space $C(\R_+,\R_+)$ of
continuous functions from $\R_+$ into $\R_+$. We write
$$\sigma(e)=\inf\{s> 0:e(s)=0\}$$
for the duration of excursion $e$. 
For $s>0$, $n_{(s)}$ will denote the conditioned measure $n(\cdot\mid\sigma=s)$.
Our normalization of the excursion measure is fixed by the relation
$$n=\int_0^\infty {ds\over 2\sqrt{2\pi s^3}}\;n_{(s)}.$$

 If $x\in\R^d$, the excursion measure
$\N_x$ of the Brownian snake from $x$ is then defined by
$$\N_x=\int_{C(\R_+,\R_+)} n(de)\;\Theta^e_{\ov x}$$
where $\ov x$ denotes the trivial element of $\W$ with lifetime $0$
and initial point $x$. Alternatively, we can view $\N_x$ as
the excursion measure of the Brownian snake from the 
regular point $\ov x$. With a slight abuse of notation
we will also write $\sigma(\omega)=\inf\{s>0:\zeta_s(\omega)=0\}$ for
$\omega\in\Omega$. We can then consider the conditioned measures
$$\N_x^{(s)}=\N_x(\cdot\mid \sigma=s)=
\int_{C(\R_+,\R_+)} n_{(s)}(de)\;\Theta^e_{\ov x}.$$
Note that in contrast to the introduction we now view $\N_x^{(s)}$
as a measure on $\Omega$ rather than on $C([0,s],\W)$.
The range $\r=\r(\omega)$ is defined by
$\r=\{\wh W_s:s\geq 0\}.$

\begin{lemma}
\label{hitting}
Suppose that $d=1$ and let $x>0$. {\rm(i)} We have
$$\N_x(\r\cap ]-\infty,0]\not =\varnothing )={3\over 2x^2}.$$
{\rm(ii)} For every $\lambda>0$,
$$\N_x\Big(1-\ind{\{\r\cap ]-\infty,0]=\varnothing \}}\,e^{-\lambda\sigma}\Big)=
\sqrt{{\lambda\over 2}}\Big(3\,({\rm
coth}(2^{1/4}x\lambda^{1/4}))^2-2\Big)$$
where ${\rm
coth}(y)=\cosh(y)/\sinh(y)$. 
\end{lemma}

\proof (i) According to Section VI.1 of \cite{Zu}, the function
$u(x)=\N_x(\r\cap ]-\infty,0]\not =\varnothing )$ solves $u''=4\, u^2$ in
$]0,\infty[$, with boundary condition $u(0+)=+\infty$. The
desired result follows.

(ii) See Lemma 7 in \cite{Delmas}. \cq

\subsection{Finite-dimensional marginal distributions}

In this subsection we state a result giving information
about the joint distribution of the values of the Brownian snake at a finite
number of times and its range. In order to state this result, we need 
some formalism for trees. We first introduce the set of labels
$${\cal U}=\bigcup_{n=0}^\infty \{1,2\}^n  $$
where by convention $\{1,2\}^0=\{\varnothing \}$. An
element
of ${\cal U}$ is thus a sequence
$u=u^1\ldots u^n$ of elements of $\{1,2\}$, and we set $|u|=n$, so that 
$|u|$ represents the ``generation'' of $u$. In particular, $|\varnothing|=0$.
The mapping $\pi:{\cal U}\backslash\{\varnothing \}\la {\cal U}$ is 
defined by $\pi(u^1\ldots u^n)=u^1\ldots u^{n-1}$ ($\pi(u)$ is the
``father'' of $u$). In particular, if $k=|u|$, we have $\pi^k(u)=\varnothing$.

A binary (plane) tree $\t$ is a finite subset of
$\cal U$ such that:
\begin{description}
\item{(i)} $\varnothing \in \t$.

\item{(ii)} $u\in \t\backslash\{\varnothing \}\Rightarrow
\pi(u)\in \t$.

\item{(iii)} For every $u\in\t$, either $u1\in \u$ and $u2\in\u$,
or $u1\notin \u$ and $u2\notin \u$ ($u$ is called a leaf in the second case).
\end{description}

We denote by ${\bf A}$ the set of all binary trees. A marked tree is then
a pair $(\t,(h_u)_{u\in\t})$ where $\t\in{\bf A}$ and $h_u\geq 0$ for every $u\in \t$.
We denote by $\T$ the space of all marked trees.
In this work it will be convenient to view marked trees as 
$\R$-trees in the sense of \cite{DuLG} 
or \cite{EPW} (see also Section 1 above). This can be achieved through the following explicit
construction. Let $\theta=(\t,(h_u)_{u\in\t})$ be a marked tree and let 
$\R^\t$ be the vector space of all mappings from $\t$ into $\R$. Write 
$(\varepsilon_u,u\in\t)$ for the canonical basis of $\R^\t$. Then consider the
mapping 
$$p_\theta:\bigcup_{u\in\t}\{u\}\times[0,h_u]\la \R^\t$$
defined by
$$p_\theta(u,\ell)=\sum_{k=1}^{|u|} h_{\pi^k(u)}\,\varepsilon_{\pi^k(u)} +\ell\,\varepsilon_u.$$
As a set, the real tree associated with $\theta$ is  the range $\wt\theta$
of $p_\theta$. Note that this is a connected union of line segments in $\R^\t$.
It is equipped with the distance $d_{\theta}$ such that $d_{\theta}(a,b)$
is the length of the shortest path in $\wt \theta$ going from $a$ to $b$. By
definition, the range of this path is the segment between $a$ and $b$ and is denoted
by
$\llbracket a,b\rrbracket$.  
Finally, we will write $\l_\theta$ for (one-dimensional) Lebesgue measure
on $\wt \theta$. 

By definition, leaves of $\wt\theta$ are points of the form
$p_\theta(u,h_u)$ where $u$ is leaf of $\theta$. Points of the form
$p_\theta(u,h_u)$ when $u$ is not a leaf are called nodes of $\wt \theta$.
We write $L(\theta)$ for the set of leaves of $\wt \theta$, and 
$I(\theta)$ for the set of its nodes. The root of
$\wt \theta$ is just the point $0=p_\theta(\varnothing ,0)$.

We will consider Brownian motion indexed by $\wt \theta$, with initial point
$x\in\R^d$.
Formally, we may consider, under the probability measure $Q^\theta_x$, a 
collection $(\xi^u)_{u\in\t}$ of independent $d$-dimensional Brownian motions
all started at $0$ except $\xi^\varnothing$
which starts at $x$, and define a continuous process
$(V_a,a\in\wt\theta)$ by setting
$$V_{p_\theta(u,\ell)}=\sum_{k=1}^{|u|} \xi^{\pi^k(u)}(h_{\pi^k(u)})+\xi^u(\ell),$$
for every $u\in\t$ and $\ell\in[0,h_u]$.
Finally, with every leaf $a$ of $\wt\theta$ we associate a stopped path
$\w^{(a)}$ with lifetime $d_\theta(0,a)$: For every $t\in[0,d_\theta(0,a)]$,
$\w^{(a)}(t)=V_{r(a,t)}$ where $r(a,t)$ is the unique element of $\llbracket
0,a\rrbracket$ such that
$d_\theta(0,r(a,t))=t$.

For every integer $p\geq 1$, denote by ${\bf A}_p$ the set of all binary trees
with $p$ leaves, and by $\T_p$ the corresponding set of marked trees. The uniform
measure $\Lambda_p$ on $\T_p$ is defined by
$$\int_{\T_p} \Lambda_p(d\theta)\,F(\theta)=
\sum_{\t\in{\bf A}_p} \int \prod_{v\in\t} dh_v\,F(T,(h_v)_{v\in\t}).$$

With this notation, Proposition IV.2 of \cite{Zu} states that, for
every integer $p\geq 1$ and every 
symmetric
nonnegative measurable function $F$ on $\W^p$,
\be
\label{margi}
\N_x\Big(\int_{]0,\sigma[^p} ds_1\ldots
ds_p\,F(W_{s_1},\ldots,W_{s_p})\Big) =2^{p-1}\,p! \int
\Lambda_p(d\theta)\,Q^\theta_x\Big[F((\w^{(a)})_{a\in L(\theta)})\Big].
\ee
We will need a stronger result concerning the case where the function $F$
also depends on the range $\r$ of the Brownian snake. To state this result, denote
by ${\cal K}$ the space of all compact subsets of $\R^d$, which is equipped
with the Hausdorff metric and the associated Borel $\sigma$-field. Suppose that
under the probability measure $Q^\theta_x$ (for each choice of $\theta$ in $\T$),
in addition to the process $(V_a,a\in\wt\theta)$,
we are also given an independent Poisson point measure on $\wt\theta\times \Omega$,
denoted by
$$\sum_{i\in I} \delta_{(a_i,\omega_i)},$$
with intensity $4\,{\cal L}_\theta(da)\otimes \N_0(d\omega)$.

\begin{theorem}
\label{marginals}
For every nonnegative measurable function $F$
on $\W^p\times {\cal K}\times\R_+$, which is symmetric
in the first $p$ variables, we have
\ba&&\N_x\Big(\int_{]0,\sigma[^p} ds_1\ldots
ds_p\,F(W_{s_1},\ldots,W_{s_p},\r,\sigma)\Big)\\ 
&&\quad=2^{p-1}\,p! \int
\Lambda_p(d\theta)\,Q^\theta_x\Big[F\Big((\w^{(a)})_{a\in L(\theta)},
{\rm cl}\Big(\bigcup_{i\in I}\; (V_{a_i}+\r(\omega_i))\Big),\sum_{i\in I} \sigma(\omega_i)\Big)\Big],
\ea
where ${\rm cl}(A)$ denotes the closure of the set $A$.
\end{theorem}

\noindent{\bf Remark.} It is immediate to see that
$${\rm cl}\Big(\bigcup_{i\in I}\; (V_{a_i}+\r(\omega_i))\Big)=
\Big(\bigcup_{a\in L(\theta)}\w^{(a)}[0,\zeta_{(\w^{(a)})}]
\Big)\cup\Big(\bigcup_{i\in I}\; (V_{a_i}+\r(\omega_i))\Big),\qquad Q^\theta_x\hbox{
a.e.}$$

\proof Consider first the case $p=1$. Let $F_1$ be a nonnegative measurable function
on
$\W$, and let $F_2$ and $F_3$ be two nonnegative measurable functions on $\Omega$.
By applying the Markov property under $\N_x$ at time $s$, then using the
time-reversal invariance of $\N_x$ (which is easy from the analogous
property for the It\^o measure $n(de)$), and finally using the Markov property at
time
$s$ once again, we get
\ba
&&\N_x\Big(\int_0^\sigma ds \,F_1(W_s)\,F_2\Big((W_{(s-r)^+})_{r\geq
0}\Big)\,F_3\Big((W_{s+r}) _{r\geq 0}\Big)\Big)\\
&&\qquad=\N_x\Big(\int_0^\sigma ds \,F_1(W_s)\,F_2\Big((W_{(s-r)^+})_{r\geq
0}\Big)\,\E_{W_s}\Big[F_3\Big((W_{r\wedge \sigma})_{r\geq 0}\Big)\Big]\Big)\\
&&\qquad=\N_x\Big(\int_0^\sigma ds \,F_1(W_s)\,F_2\Big((W_{s+r})_{r\geq
0}\Big)\,\E_{W_s}\Big[F_3\Big((W_{r\wedge \sigma})_{r\geq 0}\Big)\Big]\Big)\\
&&\qquad=\N_x\Big(\int_0^\sigma ds \,F_1(W_s)\,
\E_{W_s}\Big[F_2\Big((W_{r\wedge\sigma})_{r\geq
0}\Big)\Big]\,\E_{W_s}\Big[F_3\Big((W_{r\wedge \sigma})_{r\geq 0}\Big)\Big]\Big).
\ea
We then use the case $p=1$ of (\ref{margi}) to see that the last quantity
is equal to
$$\int_0^\infty dt\,\int P_x^t(d\w)\,F_1(\w)\,
\E_\w\Big[F_2\Big((W_{r\wedge\sigma})_{r\geq
0}\Big)\Big]\,\E_\w\Big[F_3\Big((W_{r\wedge\sigma})_{r\geq
0}\Big)\Big],$$
where $P_x^t$ denotes the law of Brownian motion started at $x$ and stopped at time
$t$ (this law is viewed as a probability measure on $\W$). Now if we specialize to
the case where $F_2$ is a function of the form $F_2(\omega)
=G_2(\{\wh W_s(\omega):s\geq 0\},\sigma)$, an immediate application of Lemma V.2 
in \cite{Zu} shows that
$$\E_\w\Big[F_2\Big((W_{r\wedge\sigma})_{r\geq
0}\Big)\Big]=E\Big[G_2\Big({\rm cl}\Big(\bigcup_{j\in J}\;
(\w(t_j)+\r(\omega_j))\Big),\sum_{j\in J}\sigma(\omega_j)\Big)\Big],$$
where $\sum_{j\in J}\delta_{(t_j,\omega_j)}$ is a Poisson point measure
on $[0,\zeta_{(\w)}]\times\Omega$ with intensity $2\;dt\,\N_0(d\omega)$.
Applying the same observation to $F_3$, we easily get the case
$p=1$ of the theorem.

The general case can be derived along similar lines by using Theorem 3 in
\cite{LG1}. Roughly speaking, the case $p=1$ amounts to combining 
Bismut's decomposition of the Brownian excursion (Lemma 1 in \cite{LG1})
with the spatial displacements of the Brownian snake. For general $p$,
the second assertion of Theorem 3 in \cite{LG1} provides the analogue of Bismut's
decomposition, which when combined with spatial displacements leads to the statement
of Theorem
\ref{marginals}. Details are left to the reader. \cq

\subsection{The re-rooting theorem}

In this subsection, we state and prove an important invariance property
of the Brownian snake under $\N_0$, which plays a major role in Section 3 below.
We first need to introduce some notation. For every $s,r\in[0,\sigma]$,
we set
$$s\oplus r=\left\{\begin{array}{ll}
s+r&\hbox{if }s+r\leq \sigma\;,\\
s+r-\sigma\quad&\hbox{if }s+r> \sigma\;.
\end{array}
\right.
$$
We also use the following convenient notation for closed intervals: If
$u,v\in \R$, $[u,v]=[v,u]=[u\wedge v,u\vee v]$.

Let $s\in[0,\sigma[$. In order to define the re-rooted snake $W^{[s]}$,
we first set
$$\zeta^{[s]}_r=\zeta_s+\zeta_{s\oplus r}-2\;\inf_{u\in [s,s\oplus r]}\zeta_u\;,$$
if $r\in[0,\sigma]$, and $\zeta^{[s]}_r=0$ if $r>\sigma$. We also want to define
the stopped paths $W^{[s]}_r$, in such a way that
$$\wh W^{[s]}_r=\wh W_{s\oplus r}-\wh W_s\;,$$
if $r\in[0,\sigma]$, and $\wh W^{[s]}_r=0$ if $r>\sigma$. To this end, we may 
notice that $\wh W^{[s]}$ satisfies the property
$$\wh W^{[s]}_r=\wh W^{[s]}_{r'}\qquad{\rm if}\
\zeta^{[s]}_r=\zeta^{[s]}_{r'}=\inf_{u\in[r,r']}\zeta^{[s]}_u$$
and so in the terminology of \cite{MM2}, $(W^{[s]}_r)_{0\leq r\leq \sigma}$
is uniquely determined as the snake whose tour is $(\zeta^{[s]}_r,\wh
W^{[s]}_r)_{0\leq r\leq \sigma}$ (see the homeomorphism theorem of \cite{MM2}).
We have the explicit formula, for $r\geq 0$,
and $0\leq t\leq \zeta^{[s]}_r$,
\be
\label{snaketour}
W^{[s]}_r(t)=\wh W^{[s]}_{\sup\{u\leq r\;:\;\zeta^{[s]}_u=t\}}.
\ee
As explained in the introduction,
$(\zeta^{[s]}_r)_{r\geq 0}$ codes the same real tree as the one coded by
$(\zeta_r)_{r\geq 0}$, but with  a new root which is the vertex originally labelled
by $s$, and $\wh W^{[s]}_r$ gives  the
spatial displacements along the line segment from the (new) root to the vertex coded 
by $r$ (in the coding given by $\zeta^{[s]}$).

\begin{theorem}
\label{reroot}
For every nonnegative measurable function $F$
on $\R_+\times \Omega$,
$$\N_0\Big(\int_0^\sigma ds\,F(s,W^{[s]})\Big)
=\N_0\Big(\int_0^\sigma ds\,F(s,W)\Big).$$
\end{theorem}

\noindent {\bf Remark.} For every $s\in[0,\sigma[$, the duration of the 
re-rooted snake excursion $W^{[s]}$ is the same as that of
the original one. Using this simple observation, and replacing 
$F$ by $\ind{\{1-\varepsilon<\sigma\leq 1\}}F$, we can easily get
a version of Theorem \ref{reroot} for the normalized
Brownian snake excursion. Precisely,
the formula of Theorem \ref{reroot} still holds if $\N_0$ is
replaced by $\N^{(1)}_0$ (or by $\N^{(r)}_0$ for any $r>0$). 
Via a continuity argument, it follows that, for every $s\in[0,1[$,
and every nonnegative measurable function $G$ on $\Omega$,
$$\N^{(1)}_0\left(G(W^{[s]})\right)=\N^{(1)}_0(G(W)).$$
See 
Proposition 13 in \cite{MM} for the same result
with a different approach.

\proof By (\ref{snaketour}), $W^{[s]}$ can be written $\N_0$ a.e. as
$\Phi(\zeta^{[s]},\wh W^{[s]})$, where the deterministic function $\Phi$ does not
depend on $s$. Also note that when $s=0$, $W=W^{[0]}=\Phi(\zeta,\wh W)$, $\N_0$ a.e.
In view of these considerations, it will be sufficient to treat the case when
$$F(s,W)=F_1(s,\zeta)\,F_2(s,\wh W)$$
where $F_1$ and $F_2$ are nonnegative measurable functions defined 
respectively on $\R_+\times C(\R_+,\R_+)$ and on $\R_+\times C(\R_+,\R^d)$. We first
deal with the special case $F_2=1$. 

For $s\in[0,\sigma[$ and $r\geq 0$, set
\ba
&&\zeta^{1,s}_r=\zeta_{(s-r)^+}-\zeta_s\;,\\
&&\zeta^{2,s}_r=\zeta_{s+r}-\zeta_s\;.
\ea
Let $G$ be a nonnegative measurable function 
on $\R_+\times C(\R_+,\R)\times\R_+\times C(\R_+,\R)$. From the Bismut decomposition
of the Brownian excursion (see e.g. Lemma 1 in \cite{LG1}), we have
$$\N_0\Big(\int_0^\sigma ds\,G\Big(s,(\zeta^{1,s}_r)_{r\geq
0},\sigma-s,(\zeta^{2,s}_r)_{r\geq 0}\Big)\Big)
=\int_0^\infty da\,E\Big[G\Big(T_a,(B_{r\wedge T_a})_{r\geq 0},T'_a,(B'_{r\wedge
T'_a})_{r\geq 0}\Big)\Big],$$
where $B$ and $B'$ are two independent linear Brownian motions
started at $0$, and
$$T_a=\inf\{r\geq 0:B_r=-a\}\ ,\quad T'_a=\inf\{r\geq 0:B'_r=-a\}.$$
Now observe that
\ba
&&\zeta^{[s]}_r=\zeta^{2,s}_r-2\,\inf_{0\leq u\leq r} \zeta^{2,s}_u\ ,\qquad
\hbox{if }0\leq r\leq \sigma-s\;,\\
&&\zeta^{[s]}_{\sigma-r}=\zeta^{1,s}_r-2\,\inf_{0\leq u\leq r} \zeta^{1,s}_u\
,\qquad
\hbox{if }0\leq r\leq s\;,
\ea
and note that $R_t:=B_t-2\inf_{r\leq t}B_r$ and $R'_r:=B'_t-2 \inf_{r\leq t} B'_r$
are two independent three-dimensional Bessel processes, for which
\ba &&L_a:=\sup\{t\geq 0:R_t\leq a\}=T_a\;,\\
&&L'_a:=\sup\{t\geq 0:R'_t\leq a\}=T'_a\;.
\ea
(This is Pitman's theorem, see e.g. \cite{RY}, Theorem VI.3.5.) It follows that
\ba
&&\N_0\Big(\int_0^\sigma ds\,G\Big(\sigma-s,(\zeta^{[s]}_{r\wedge(\sigma-s)})_{r\geq
0},s,(\zeta^{[s]}_{\sigma-(r\wedge s)})_{r\geq 0}\Big)\Big)\\
&&\qquad = \int_0^\infty da \,E\Big[G\Big(L'_a,(R'_{r\wedge L'_a})_{r\geq 0},
L_a,(R_{r\wedge L_a})_{r\geq 0}\Big)\Big)\\
&&\qquad=\N_0\Big(\int_0^\sigma ds\,G\Big(s,(\zeta_{r\wedge s})_{r\geq 0},\sigma-s,
(\zeta_{(\sigma-r)\vee s})_{r\geq 0}\Big)\Big)
\ea
where the last equality is again a consequence of the Bismut decomposition, together
with the Williams reversal theorem (\cite{RY}, Corollary XII.4.4). 
Changing $s$ into $\sigma-s$ in the last
integral gives the desired result when $F_2=1$. 

Let us consider the general case. For simplicity we take $d=1$, but the
argument can obviously be extended. From the definition of the Brownian snake,
we have
$$\N_0\Big(\int_0^\sigma ds\,F_1(s,\zeta)\,F_2(s,\wh W)\Big)
=\N_0\Big(\int_0^\sigma ds\,F_1(s,\zeta)\,\Theta^\zeta_0[F_2(s,\wh W)]\Big)$$
and $\wh W$ is under $\Theta^\zeta_0$ a centered Gaussian process with covariance
$${\rm cov}_{\Theta^\zeta_0}(\wh W_s,\wh W_{s'})=\inf_{r\in[s,s']}\zeta_r.$$
We have in particular
$$\N_0\Big(\int_0^\sigma ds\,F_1(s,\zeta^{[s]})\,F_2(s,\wh W^{[s]})\Big)
=\N_0\Big(\int_0^\sigma ds\,F_1(s,\zeta^{[s]})\,\Theta^\zeta_0\Big[F_2(s,
(\wh W_{s\oplus r}-\wh W_s)_{r\geq 0})\Big]\Big).$$
Now note that $(\wh W_{s\oplus r}-\wh W_s)_{r\geq 0}$ is under $\Theta^\zeta_0$
a Gaussian process with covariance
$${\rm cov}(\wh W_{s\oplus r}-\wh W_s,\wh W_{s\oplus r'}-\wh W_s)
=\inf_{[s\oplus r,s\oplus r']}\zeta_u-\inf_{[s\oplus r,s]}\zeta_u-
\inf_{[s\oplus r',s]}\zeta_u+\zeta_s=\inf_{[r,r']}\zeta^{[s]}_u,$$
where the last equality follows from an elementary verification. Hence,
$$\Theta^\zeta_0\Big[F_2(s,
(\wh W_{s\oplus r}-\wh W_s)_{r\geq 0})\Big]=
\Theta^{\zeta^{[s]}}_0[F_2(s,\wh W)],$$
and, using the first part of the proof,
\ba
\N_0\Big(\int_0^\sigma ds\,F_1(s,\zeta^{[s]})\,F_2(s,\wh W^{[s]})\Big)&=&
\N_0\Big(\int_0^\sigma ds\,F_1(s,\zeta^{[s]})\,
\Theta^{\zeta^{[s]}}_0[F_2(s,\wh W)]\\
&=&\N_0\Big(\int_0^\sigma ds\,F_1(s,\zeta)\,
\Theta^{\zeta}_0[F_2(s,\wh W)]\Big)\\
&=&\N_0\Big(\int_0^\sigma ds\,F_1(s,\zeta)\,F_2(s,\wh W)\Big).
\ea
This completes the proof. \cq

\subsection{The special Markov property}

Let $D$ be a domain in $\R^d$, and fix a point $x\in D$. For every
$\w\in\W$, we set
$$\tau(\w):=\inf\{t\geq 0:\w(t)\notin D\}$$
where $\inf\varnothing =+\infty$ as usual. The random set
$$\{s\geq 0:\tau(W_s)<\zeta_s\}$$
is open $\N_x$ a.e., and can thus be written as a disjoint 
union of open intervals $]a_i,b_i[$, ${i\in I}$.
It is easy to verify that $\N_x$ a.e. for every $i\in I$
and every $s\in]a_i,b_i[$,
$$\tau(W_s)=\tau(W_{a_i})=\tau(W_{b_i})=\zeta_{a_i}=\zeta_{b_i}$$
and moreover the paths $W_s,\,s\in[a_i,b_i]$ coincide up to their
exit time from $D$.

For every $i\in I$, we define a random element $W^{(i)}$ of $\Omega$
by setting for every $s\geq 0$
$$W^{(i)}_s(t)=W_{(a_i+s)\wedge b_i}(\zeta_{a_i}+t)\,,\quad \hbox{for }
0\leq t\leq \zeta_{(W^{(i)}_s)}:=\zeta_{(a_i+s)\wedge b_i}-\zeta_{a_i}.$$
Informally, the $W^{(i)}$'s represent the excursions of the 
Brownian snake outside $D$ (the word ``outside'' is a bit misleading
since these excursions may come back into $D$ even though they start
from the boundary of $D$).

Finally, we also need a process that contains the information given by
the Brownian snake paths before they exit $D$. We set
$\wt W_s^D=W_{\eta^D_s}$, where for every $s\geq 0$,
$$\eta^D_s:=\inf\{r\geq 0:\int_0^r du\,\ind{\{\tau(W_u)\geq \zeta_u\}}>s\}.$$
The $\sigma$-field $\e^D$ is by definition generated by the process $\wt W^D$
and by the class of $\N_x$-negligible subsets of $\Omega$
(the point $x$ is fixed throughout this subsection). The following 
statement is proved
in \cite{LG2} (Proposition 2.3 and Theorem 2.4).

\begin{theorem}
\label{SMP}
There exists a random finite measure denoted by $\z^D$, which is
$\e^D$-measurable and $\N_x$ a.e. supported on $\partial D$, such that
the following holds. Under $\N_x$, conditionally on $\e^D$, the point measure
$$\n:=\sum_{i\in I} \delta_{(W^{(i)})}$$
is Poisson with intensity $\int_{\partial D} \z^D(dy)\,\N_y(\cdot)$.
\end{theorem}

We will apply this theorem to the case $d=1$, $x=0$ and $D=]c,\infty[$
for some $c<0$. In that case, the measure $\z^D$ is a random multiple of the
Dirac measure at $c$: $\z^D=L^c\,\delta_c$ for some nonnegative
random variable $L^c$. From Lemma \ref{hitting}(i) and Theorem \ref{SMP}, it is easy
to verify that
$\{L^c>0\}=\{\r\cap]-\infty,c]
\not =\varnothing \}=\{\r\cap]-\infty,c[
\not =\varnothing \}\;$, $\N_0$ a.e. Moreover, as a simple consequence of the
special Markov property, the process $(L^{-r})_{r>0}$ is a nonnegative 
martingale under $\N_0$ (it is indeed a critical
continuous-state branching process). In particular the variable
$$L^{*,r}:=\sup_{c\in\Q\cap]-\infty,r]} L^c$$
is finite $\N_0$ a.e., for every $r<0$. 

\subsection{Uniqueness of the minimum}

From now on we assume that $d=1$. In this subsection, we consider the Brownian snake
under  its excursion measure $\N_0$. We use the notation
$$\un W=\inf_{s\geq 0} \wh W_s.$$
Note that the law of $\un W$ under $\N_0$ is given by Lemma \ref{hitting}(i)
and an obvious translation argument.

\begin{proposition}
\label{unimin}
There exists $\N_0$ a.e. a unique instant $s_*\in]0,\sigma[$ such that
$\wh W_{s_*}=\un W$.
\end{proposition}

This result already appears as Lemma 16 in \cite{MM}, where its proof is attributed
to T. Duquesne. We provide a short proof for the sake of completeness and also
because this result plays a major role throughout this work.

\proof Set
$$\lambda:=\inf\{s\geq 0:\wh W_s=\un W\}\ ,\quad \rho=\sup\{s\geq 0:\wh W_s=\un W\}$$
so that $0<\lambda\leq \rho<\sigma$. We have to prove that $\lambda=\rho$. 
To this end we fix $\delta>0$ and  we verify that $\N_0(\rho-\lambda>\delta)=0$.

Fix two rational numbers $q<0$ and $\varepsilon>0$. We first get an upper bound
on the quantity
$$\N_0(q-\varepsilon\leq \un W<q\,,\,\rho-\lambda>\delta).$$
Denote by $(W^{(i)})_{i\in I}$ the excursions of the Brownian snake outside
$]q,\infty[$, and by $\n$
the corresponding point measure, as in the previous subsection. Since the law of
$\un W$ under $\N_0$ has no atoms, the numbers $\un W^{(i)}$, $i\in I$
are distinct $\N_0$ a.e. Therefore, on the event $\{\un W<q\}$, the
whole interval $[\rho,\lambda]$ must be contained 
in a single excursion interval below level $q$. Hence,
$$\N_0(q-\varepsilon\leq \un W<q\,,\,\rho-\lambda>\delta)
\leq \N_0(\{\forall i\in I:\un W^{(i)}\geq q-\varepsilon\}
\cap\{\exists i\in I: \sigma(W^{(i)})>\delta\}).$$
Introduce the events $A_\varepsilon:=\{\un W<-q-\varepsilon\}$ and 
$B_{\varepsilon,\delta}
=\{\un W\geq -q-\varepsilon\;,\;\sigma>\delta\}$. We get
$$\N_0(q-\varepsilon\leq \un W<q\,,\,\rho-\lambda>\delta)
\leq \N_0(\n(A_\varepsilon)=0\,,\;\n(B_{\varepsilon,\delta})\geq 1).$$
From the special Markov property (and the remarks of the end of subsection 2.3),
we know that conditionally on $L^q$, $\n$ is a Poisson point measure
with intensity $L^q\,\N_q$. Since the sets $A_\varepsilon$ and
$B_{\varepsilon,\delta}$ are disjoint, independence properties of
Poisson point measures give
\ba
\N_0(q-\varepsilon\leq \un W<q\,,\,\rho-\lambda>\delta)
&\leq&\N_0\Big(\ind{\{\n(A_\varepsilon)=0\}}\,
(1-\exp(-L^q\N_q(B_{\varepsilon,\delta})))\Big)\\
&=&\N_0\Big(\ind{\{q-\varepsilon\leq \un W<q\}}\,
(1-\exp(-L^q\N_q(B_{\varepsilon,\delta})))\Big)\\
&\leq&\N_0\Big(\ind{\{q-\varepsilon\leq \un W<q\}}\,
(1-\exp(-c(\varepsilon,\delta)\,L^{*,q}))\Big)
\ea
where $c(\varepsilon,\delta)=\N_q(B_{\varepsilon,\delta})$ does not depend on 
$q$ by an obvious translation argument.

We can apply the preceding bound with $q$ replaced
by $q-\varepsilon,q-2\varepsilon,$ etc. By summing the resulting estimates we get
$$\N_0(\un W<q\,,\,\rho-\lambda>\delta)\leq 
\N_0\Big(\ind{\{\un W<q\}}\,
(1-\exp(-c(\varepsilon,\delta)\,L^{*,q}))\Big).$$
Clearly $c(\varepsilon,\delta)$ tends to $0$ as $\varepsilon\to 0$, and 
dominated convergence gives $\N_0(\un W<q\,,\,\rho-\lambda>\delta)=0$. This
completes the proof since $q$ was arbitrary. \cq

\subsection{Bessel processes}

Throughout this work, $(\xi_t)_{t\geq 0}$ will stand for a linear
Brownian motion started at $x$ under the probability measure $P_x$.
The notation $\xi[0,t]$ will stand for the range of $\xi$
over the time interval $[0,t]$.
For every $\delta>0$, $(R_t)_{t\geq 0}$ will denote a Bessel process 
of dimension $\delta$ started at $x$ under the probability measure $P^{(\delta)}_x$.
We will use repeatedly the following simple facts. First, if $\lambda>0$, the process
$R^{(\lambda)}_t:=\lambda^{-1}R_{\lambda^2 t}$ is under $P^{(\delta)}_x$ is Bessel process
of dimension $\delta$ started at $x/\lambda$. Secondly, if $0\leq x\leq x'$ and $t\geq 0$, the law 
of $R_t$ under $P^{(\delta)}_x$ is stochastically bounded by the law of $R_t$ under $P^{(\delta)}_{x'}$.
The latter fact follows from standard comparaison theorems applied to squared Bessel processes.

Absolute continuity relations between Bessel processes, which are consequences of
the Girsanov theorem, were first observed by Yor \cite{Yor}. We state a special case of these
relations, which will play an important role in this work. This special case appears
in Exercise XI.1.22 of \cite{RY}.

\begin{proposition}
\label{Bessel}
Let $t>0$ and let $F$ be a nonnegative measurable function on $C([0,t],\R)$. 
Then, for every $x>0$ and $\lambda>0$,
$$E_x\Big[\ind{\{\xi[0,t]\subset]0,\infty[\}}\,
\exp\Big(-{\lambda^2\over 2}\int_0^t {dr\over \xi_r^2}\Big)
\,F((\xi_r)_{0\leq r\leq t})\Big]
=x^{\nu+{1\over 2}}\,E_x^{(2+2\nu)}\Big[(R_t)^{-\nu-{1\over 2}}\,F((R_r)_{0\leq r\leq
t})\Big],$$
where $\nu=\sqrt{\lambda^2+{1\over 4}}$.
\end{proposition}

We shall be concerned by the case when $\lambda^2/2=6$, and then $2+2\nu=9$ and $\nu+{1/2}=4$.
Taking $F=1$ in that case, we see that 
\be
\label{easybound}
x^{4}\,E^{(9)}_x[R_t^{-4}]\leq 1.
\ee

\section{Conditioning and re-rooting of trees}

This section contains the proof of 
Theorem \ref{main-cond} and Theorem \ref{verwaat}
which were stated in the introduction. Both will be derived 
as consequences of Theorem \ref{exc-cond} below. Recall the notation
$s_*$ for the unique time of the minimum of $\wh W$ under $\N_0$,
and $W^{[s]}$ for the snake re-rooted at $s$.

\begin{theorem}
\label{exc-cond}
Let $\varphi:\R_+\la\R_+$ be a continuous function such that
$\varphi(s)\leq C(1\wedge s)$ for some finite constant $C$. Let
$F:\Omega\la \R_+$ be a bounded continuous
function. Then,
$$\lim_{\varepsilon\to 0}\varepsilon^{-4}\;
\N_0\Big(\sigma\,\varphi(\sigma)\,F(W)\,
\ind{\{\un W>-\varepsilon\}}\Big)={2\over 21}\;
\N_0(\varphi(\sigma)\,F(W^{[s_*]})).$$
\end{theorem}

The proof of Theorem \ref{exc-cond} occupies most of the remainder
of this section. This proof will depend on a series of lemmas. 
To motivate these lemmas, we first observe that, from the 
re-rooting identity Theorem \ref{reroot}, we have
\begin{eqnarray}
\label{reroot1}
\N_0\Big(\sigma\,\varphi(\sigma)\,F(W)\,
\ind{\{\un W>-\varepsilon\}}\Big)&=&
\N_0\Big(\varphi(\sigma)\int_0^\sigma ds\,F(W^{[s]})\,
\ind{\{\un W^{[s]}>-\varepsilon\}}\Big)\nonumber\\
&=&\N_0\Big(\varphi(\sigma)\int_0^\sigma ds\,F(W^{[s]})\,
\ind{\{\wh W_s<\un W+\varepsilon\}}\Big),
\end{eqnarray}
since $\un W^{[s]}=\un W-\wh W_s$ by construction. The fact that the 
minimum of $\wh W$ is attained at a unique time implies that
\be
\label{unimin1}
\sup\{|s-s_*|:\wh W_s<\un W+\varepsilon\} \build{\la}_{\varepsilon\to 0}
^{} 0\;,\qquad \N_0\hbox{ a.e.}
\ee
and it easily follows that
\be
\label{exccond1}
\sup\{|F(W^{[s]})-F(W^{[s_*]})|:\wh W_s
<\un W+\varepsilon\}\build{\la}_{\varepsilon\to 0}
^{} 0\;,\qquad \N_0\hbox{ a.e.}
\ee
Coming back to (\ref{reroot1}), this suggests to study the behavior
of 
$$\int_0^\sigma ds\,\ind{\{\wh W_s<\un W+\varepsilon\}}$$
as $\varepsilon\to 0$. This motivates the next two lemmas.

\begin{lemma}
\label{exccondlem1}
We have
$$\lim_{\delta\da 0}
\sup_{0<\varepsilon<1} \left(\varepsilon^{-4}\,
\N_0\Big((1-e^{-\sigma})\int_0^\sigma ds\,
\ind{\{\wh W_s<\un W+\varepsilon,\;\zeta_s<\delta\}}\Big)\right)=0.$$
\end{lemma}

\proof From Theorem \ref{reroot} and the fact that
$\zeta_s=\zeta_{\sigma-s}^{[s]}$ for every $s\in[0,\sigma]$, we have
\begin{eqnarray}
\label{rerootech8}
\N_0\Big((1-e^{-\sigma})\int_0^\sigma ds\,
\ind{\{\wh W_s<\un W+\varepsilon,\;\zeta_s<\delta\}}\Big)
&=&\N_0\Big((1-e^{-\sigma})\int_0^\sigma ds\,
\ind{\{\un W^{[s]}>-\varepsilon,\;\zeta^{[s]}_{\sigma-s}<\delta\}}\Big)
\nonumber\\
&=&\N_0\Big((1-e^{-\sigma})\int_0^\sigma ds\,
\ind{\{\un W>-\varepsilon,\;\zeta_{\sigma-s}<\delta\}}\Big)
\nonumber\\
&=&\N_0\Big((1-e^{-\sigma})\ind{\{\un W>-\varepsilon\}}\int_0^\sigma ds\,
\ind{\{\zeta_{s}<\delta\}}\Big).
\end{eqnarray}

Recall that $(\xi_t)_{t\geq 0}$ denotes a standard linear Brownian motion 
that starts from $x$ under the probability measure $P_x$, and
write $\un\xi_t:=\inf\{\xi_r:r\leq t\}$. From the case $p=1$
of Theorem \ref{marginals} and Lemma \ref{hitting}(i)
we get
\ba
\N_0\Big(\ind{\{\un W>-\varepsilon\}}\int_0^\sigma ds\,
\ind{\{\zeta_{s}<\delta\}}\Big)
&=&\int_0^\delta da\,E_0\Big[\ind{\{\un \xi_a>-\varepsilon\}}
\exp\Big(-4\int_0^a dt\,\N_{\xi_t}(\un W_t\leq -\varepsilon)\Big)\Big]\\
&=&\int_0^\delta da\,E_\varepsilon\Big[\ind{\{\un \xi_a>0\}}
\exp\Big(-6\int_0^a {dt\over \xi_t^2}\Big)\Big].
\ea
At this point we use Proposition \ref{Bessel}, which gives
$$E_\varepsilon\Big[\ind{\{\un \xi_a>0\}}
\exp\Big(-6\int_0^a {dt\over \xi_t^2}\Big)\Big]=\varepsilon^4
E^{(9)}_\varepsilon[R_a^{-4}].$$

Similarly,
\ba
\N_0\Big(\ind{\{\un W>-\varepsilon\}}e^{-\sigma}\int_0^\sigma ds\,
\ind{\{\zeta_{s}<\delta\}}\Big)
&=&\int_0^\delta da\,E_0\Big[\ind{\{\un \xi_a>-\varepsilon\}}
\exp\Big(-4\int_0^a dt\,\N_{\xi_t}(
1-\ind{\{\un W> -\varepsilon)\}}e^{-\sigma})\Big)\Big]\\
&=&\int_0^\delta da\,E_\varepsilon\Big[\ind{\{\un \xi_a>0\}}
\exp\Big(-4\int_0^a dt\,\N_{\xi_t}(
1-\ind{\{\un W> 0)\}}e^{-\sigma})\Big)\Big]\\
&=&\int_0^\delta da\,E_\varepsilon\Big[\ind{\{\un \xi_a>0\}}
\exp\Big(-2^{3/2}\int_0^a dt\,(3{\rm coth}(2^{1/4}\xi_t)^2-2)
\Big)\Big]
\ea
using Lemma \ref{hitting}(ii) in the last equality. For every $x>0$, set
$$h(x)=-{3\over 2x^2}+2^{-1/2}(3{\rm coth}(2^{1/4}x)^2-2)>0.$$
A Taylor expansion shows that $h(x)\leq C\,x^2$. (Here and later,
$C,C',C''$ denote constants whose exact value is unimportant.) Then,
\ba
\N_0\Big(\ind{\{\un W>-\varepsilon\}}e^{-\sigma}\int_0^\sigma ds\,
\ind{\{\zeta_{s}<\delta\}}\Big)
&=&\int_0^\delta da\,E_\varepsilon\Big[\ind{\{\un \xi_a>0\}}
\exp\Big(-6\int_0^a {dt\over \xi_t^2}
-4\int_0^a dt\,h(\xi_t)\Big)\Big]\\
&=&\int_0^\delta da\;
\varepsilon^4
E^{(9)}_\varepsilon\Big[R_a^{-4}\exp\Big(-4\int_0^a dt\,h(R_t)\Big)\Big]
\ea
using Proposition \ref{Bessel} as above.

By combining the preceding calculations, we arrive at
\begin{eqnarray}
\label{lem1tech0}
\varepsilon^{-4}\N_0\Big(\ind{\{\un
W>-\varepsilon\}}(1-e^{-\sigma})\int_0^\sigma ds\,
\ind{\{\zeta_{s}<\delta\}}\Big)&=&
\int_0^\delta da\,
E^{(9)}_\varepsilon\Big[R_a^{-4}\Big(1-\exp\Big(-4\int_0^a dt\,h(R_t)
\Big)\Big)\Big]\nonumber\\
&\leq&4C\int_0^\delta da\int_0^a dt E^{(9)}_\varepsilon[R_a^{-4}R_t^2].
\end{eqnarray}

If $a\leq \delta<1/2$ and $0<t\leq a^2$, we can bound
\be
\label{lem1tech1}
E^{(9)}_\varepsilon[R_a^{-4}R_t^2]
=E^{(9)}_\varepsilon\left[R_t^2\,E^{(9)}_{R_t}[R_{a-t}^{-4}]\right]
\leq E^{(9)}_\varepsilon[R_t^2]\,E^{(9)}_{0}[R_{a-t}^{-4}]
\leq C'\,a^{-2}\;.
\ee

If $a^2<t\leq a$ we use a different argument: From the bound (\ref{easybound}), we get
\be 
\label{lem1tech2}
E^{(9)}_\varepsilon[R_a^{-4}R_t^2]
=E^{(9)}_\varepsilon\left[R_t^2\,E^{(9)}_{R_t}[R_{a-t}^{-4}]\right]
\leq E^{(9)}_\varepsilon[R_t^{-2}]\leq E^{(9)}_0[R_t^{-2}]=C''t^{-1}.
\ee

By substituting the bounds (\ref{lem1tech1}) and (\ref{lem1tech2})
into (\ref{lem1tech0}), we arrive at
$$
\varepsilon^{-4}\N_0\Big(\ind{\{\un
W>-\varepsilon\}}(1-e^{-\sigma})\int_0^\sigma ds\,
\ind{\{\zeta_{s}<\delta\}}\Big)\leq
4C\int_0^\delta da\Big(C'+C''\int_{a^2}^a t^{-1}\,dt\Big),
$$
which tends to $0$ as $\delta\to 0$. Recalling  
({\ref{rerootech8}) we see that the proof of Lemma 
\ref{exccondlem1} is complete. \cq

\begin{lemma}
\label{exccondlem2}
For every $\delta>0$, 
$$\sup_{0<\varepsilon<1}\N_0\Big(\Big(\varepsilon^{-4}\int_0^\sigma ds\,
\ind{\{\wh W_s<\un W+\varepsilon,\;\zeta_s\geq\delta\}}\Big)^2\Big)
<\infty.$$
\end{lemma}

\proof We now use the case $p=2$ of Theorem \ref{marginals} to write
\be
\label{lem2tech0}
\N_0\Big(\Big(\varepsilon^{-4}\int_0^\sigma ds\,
\ind{\{\wh W_s<\un W+\varepsilon,\;\zeta_s\geq\delta\}}\Big)^2\Big)
=4\int_{\R_+^3} da\,db\,dc\,
\ind{\{a+b\geq \delta,\;a+c\geq \delta\}}\,I^{a,b,c}_\varepsilon,
\ee
where
\ba &&I^{a,b,c}_\varepsilon
=E^{a,b,c}\Big[\ind{\{{\cal G}\subset]\gamma-\varepsilon,\infty[\}}\\
&&\qquad\exp\Big(-4\Big(\int_0^a dt\,\N_{\xi_t}(\un W>\gamma-\varepsilon)
+\int_0^b dt\,\N_{\xi'_t}(\un W>\gamma-\varepsilon)
+\int_0^c dt\,\N_{\xi''_t}(\un W>\gamma-\varepsilon)\Big)\Big)\Big],
\ea
and, under the probability measure $P^{a,b,c}$:
\begin{description}
\item{$\bullet$} $(\xi_t)_{0\leq t\leq a}$ is a linear Brownian motion
started at $0$;
\item{$\bullet$} conditionally given $(\xi_t)_{0\leq t\leq a}$,
$(\xi'_t)_{0\leq t\leq b}$ and $(\xi''_t)_{0\leq t\leq c}$ are
independent linear Brownian motions started at $\xi_a$;
\item{$\bullet$} $\gamma=\xi'_b\vee \xi''_c$;
\item{$\bullet$} ${\cal G}=\{\xi_t,0\leq t\leq a\}\cup
\{\xi'_t,0\leq t\leq b\}\cup\{\xi''_t,0\leq t\leq c\}$.
\end{description}
By Lemma \ref{hitting}(i),
$$I^{a,b,c}_\varepsilon
=E^{a,b,c}\Big[\ind{\{{\cal G}\subset]\gamma-\varepsilon,\infty[\}}
\exp\Big(-6\Big(\int_0^a {dt\over(\xi_t-\gamma+\varepsilon)^2}
+\int_0^b {dt\over(\xi'_t-\gamma+\varepsilon)^2}
+\int_0^c {dt\over(\xi''_t-\gamma+\varepsilon)^2}
\Big)\Big)\Big].$$
On the event $\{{\cal G}\subset]\gamma-\varepsilon,\infty[\}$, we have
$|\xi'_b-\xi''_c|<\varepsilon$, and
$$\begin{array}{ll}
0\leq \xi'_t-\gamma+\varepsilon\leq \xi'_t-\xi'_b+\varepsilon,\qquad
&\forall t\in[0,b],\\
\noalign{\smallskip}
0\leq \xi''_t-\gamma+\varepsilon\leq \xi''_t-\xi''_c+\varepsilon,\qquad
&\forall t\in[0,c].
\end{array}
$$
We use this to derive a first bound on $I^{a,b,c}_\varepsilon$. To write
this bound in a convenient way, we introduce the following notation:
$$\begin{array}{ll}
\eta_t=\xi_{a-t}-\xi_a,\qquad&t\in[0,a],\\
\noalign{\smallskip}
\eta'_t=\xi'_{b-t}-\xi'_b,\qquad&t\in[0,b],\\
\noalign{\smallskip}
\eta''_t=\xi''_{c-t}-\xi''_c,\qquad&t\in[0,c],
\end{array}
$$
in such a way that $\eta,\eta',\eta''$ are three independent linear
Brownian motions started at $0$ under $P^{a,b,c}$, and 
$\gamma=-\eta_a-(\eta'_b\wedge\eta''_c)$. Using this notation and the 
preceding bounds on the event $\{{\cal
G}\subset]\gamma-\varepsilon,\infty[\}$, we get
$I^{a,b,c}_\varepsilon\leq J^{a,b,c}_\varepsilon$, where
\ba
J^{a,b,c}_\varepsilon\!\!\!&=&\!\!\!
E\Big[{\bf 1}{\{|\eta'_b-\eta''_c|<\varepsilon,\eta[0,a]\subset
]-(\eta'_b\wedge\eta''_c)-\varepsilon,\infty[,
\eta'[0,b]\subset]-\varepsilon,\infty[,
\eta''[0,c]\subset]-\varepsilon,\infty[\}}\\
&&\qquad\qquad\times 
\exp\Big(-6\Big(\int_0^a {dt\over(\eta_t+(\eta'_b\wedge\eta''_c)+
\varepsilon)^2}+\int_0^b {dt\over (\eta'_t+\varepsilon)^2}
+\int_0^c {dt\over (\eta''_t+\varepsilon)^2}
\Big)\Big)\Big].
\ea
To simplify notation, we have written $E$ instead of $E^{a,b,c}$, and
$\eta[0,a]$ obviously denotes the set $\{\eta_t:0\leq t\leq a\}$,
with a similar notation for $\eta'[0,b]$ and $\eta''[0,c]$.

In the preceding formula for $J^{a,b,c}_\varepsilon$, conditioning with
respect to the pair $(\eta',\eta'')$ leads to a quantity depending
on $y=(\eta'_b\wedge\eta''_c)+\varepsilon$, of the form
\be
\label{lem2tech1}
E\Big[\ind{\{\eta[0,a]\subset]-y,\infty[\}}\,\exp\Big(-6\int_0^a
{dt\over (\eta_t+y)^2}\Big)\Big]
=E_y\Big[\ind{\{\xi[0,a]\subset]0,\infty[\}}\,\exp\Big(-6\int_0^a
{dt\over \xi_t^2}\Big)\Big]
=y^4\,E^{(9)}_y[R_a^{-4}]
\ee
using Proposition \ref{Bessel} as in the proof of Lemma \ref{exccondlem1}
above. Hence,
\ba
J^{a,b,c}_\varepsilon&=&
E\Big[{\bf 1}{\{|\eta'_b-\eta''_c|<\varepsilon,
\eta'[0,b]\subset]-\varepsilon,\infty[,
\eta''[0,c]\subset]-\varepsilon,\infty[\}}\;\\
&&\ \times\,((\eta'_b\wedge\eta''_c)+\varepsilon)^4\,
E^{(9)}_{(\eta'_b\wedge\eta''_c)+\varepsilon}[R_a^{-4}]
\exp\Big(-6\Big(\int_0^b {dt\over (\eta'_t+\varepsilon)^2}
+\int_0^c {dt\over (\eta''_t+\varepsilon)^2}
\Big)\Big)\Big].
\ea
Recall that our goal is to bound $\int da\,db\,dc\,
\ind{\{a+b\geq \delta,a+c\geq \delta\}}\,J^{a,b,c}_\varepsilon$. 
First consider the integral over the set $\{a<\delta/2\}$.
Then plainly we have $b>\delta/2$ and $c>\delta/2$, and we can use 
(\ref{easybound}) to bound
\ba
&&\int_{\{a<\delta/2\}} da\,db\,dc\,
\ind{\{a+b\geq \delta,\;a+c\geq \delta\}}\,J^{a,b,c}_\varepsilon\\
&&\ \leq {\delta\over 2}\int_{]\delta/2,\infty[^2}
db\,dc\,E\Big[\ind{\{
\eta'[0,b]\subset]-\varepsilon,\infty[,\;
\eta''[0,c]\subset]-\varepsilon,\infty[\}}
\exp\Big(-6\Big(\int_0^b {dt\over (\eta'_t+\varepsilon)^2}
+\int_0^c {dt\over (\eta''_t+\varepsilon)^2}
\Big)\Big)\Big]\\
&&\ ={\delta\over 2}\int_{]\delta/2,\infty[^2}
db\,dc\,(\varepsilon^4E^{(9)}_\varepsilon[R_b^{-4}])\,
(\varepsilon^4E^{(9)}_\varepsilon[R_c^{-4}])\\
&&\ \leq C_\delta\,\varepsilon^8.
\ea
In the last inequality, we used the fact that, for every $y>0$,
\be
\label{lem2tech2}
\int_{\delta/2}^\infty db\,E^{(9)}_y[R_b^{-4}]
\leq \int_{\delta/2}^\infty db\,E^{(9)}_0[R_b^{-4}]=C'_\delta<\infty.
\ee

We still have to get a similar bound for the integral over the
set $\{a\geq \delta/2\}$. Applying the bound (\ref{lem2tech2})
with $y=(\eta'_b\wedge\eta''_c)+\varepsilon$, we see that it is 
enough to prove that
\begin{eqnarray}
\label{lem2tech3}
&&\int_{[0,\infty[^2} db\,dc\,
E\Big[{\bf 1}{\{|\eta'_b-\eta''_c|<\varepsilon,\;
\eta'[0,b]\subset]-\varepsilon,\infty[,\;
\eta''[0,c]\subset]-\varepsilon,\infty[\}}\nonumber\\
&&\qquad\qquad\qquad\times((\eta'_b\wedge\eta''_c)+\varepsilon)^4
\exp\Big(-6\Big(\int_0^b {dt\over (\eta'_t+\varepsilon)^2}
+\int_0^c {dt\over (\eta''_t+\varepsilon)^2}
\Big)\Big)\Big]\leq C\,\varepsilon^8.
\end{eqnarray}
From Proposition \ref{Bessel} again, the left-hand side of
(\ref{lem2tech3}) is equal to
\be
\label{lem2tech4}
\varepsilon^8\int_{[0,\infty[^2} db\,dc\,
E^{(9)}_\varepsilon\otimes E^{(9)}_\varepsilon[
\ind{\{|R_b-\wt R_c|<\varepsilon\}}(R_b\wedge \wt R_c)^4
\,R_b^{-4}\,\wt R_c^{-4}],
\ee
where $R$ and $\wt R$ are two independent nine-dimensional Bessel
processes started at $\varepsilon$ under the probability measure 
$P^{(9)}_\varepsilon\otimes P^{(9)}_\varepsilon$. The quantity
(\ref{lem2tech4}) is bounded above by
$\varepsilon^8(I^\varepsilon_1+I^\varepsilon_2)$, where
\ba I^\varepsilon_1
&=&\int db\,dc\,E^{(9)}_\varepsilon\otimes E^{(9)}_\varepsilon[
\ind{\{R_b<4\varepsilon,\;\wt R_c<4\varepsilon\}}(R_b\wedge \wt R_c)^4
\,R_b^{-4}\,\wt R_c^{-4}]\\
&\leq& \Big(E^{(9)}_\varepsilon\Big[\int_0^\infty 
db\,R_b^{-2}\,\ind{\{R_b<4\varepsilon\}}
\Big]\Big)^2\\
&=&C<\infty
\ea
and 
$$I^\varepsilon_2=\int db\,dc\,
E^{(9)}_\varepsilon\otimes E^{(9)}_\varepsilon[
\ind{\{\wt R_c>3\varepsilon,\;|R_b-\wt R_c|<\varepsilon\}}(R_b\wedge \wt
R_c)^4
\,R_b^{-4}\,\wt R_c^{-4}]$$
To bound $I^\varepsilon_2$, note that, if $y\geq 3\varepsilon$,
$$E^{(9)}_\varepsilon\Big[\int_0^\infty db\,R_b^{-2}\,
\ind{\{|R_b-y|<\varepsilon\}}\Big]
=C\int_{\R^9}dz\,|z-z_\varepsilon|^{-9}\,\ind{\{y-\varepsilon
<|z|<y+\varepsilon\}}\leq C'\,{\varepsilon\over y},$$
where the notation $z_\varepsilon$ stands for a point
in $\R^9$ such that $|z_\varepsilon|=\varepsilon$, and we used the
form of the Green function of nine-dimensional Brownian motion. 
It follows that
$$I^\varepsilon_2\leq 
C'\varepsilon\,E^{(9)}_\varepsilon\Big[\int_0^\infty
dc\,R_c^{-3}\,\ind{\{R_c\geq 3\varepsilon\}}\Big]
=C'\varepsilon\,E^{(9)}_0\Big[\int_0^\infty
dc\,R_c^{-3}\,\ind{\{R_c\geq 3\varepsilon\}}\Big]
=C''<\infty,$$
by a simple scaling argument. This completes the proof 
of the bound (\ref{lem2tech3}) and of Lemma \ref{exccondlem2}. \cq

\begin{corollary}
\label{exccondcoro}
If $F:\Omega\la \R_+$ is bounded and continuous, we have
$$\lim_{\varepsilon\to 0} \varepsilon^{-4}\,
\N_0\Big((1-e^{-\sigma})\int_0^\sigma ds\,
|F(W^{[s]})-F(W^{[s_*]})|\,
\ind{\{\wh W_s<\un W+\varepsilon\}}\Big)=0.$$
\end{corollary}

\proof Thanks to Lemma \ref{exccondlem1}, it is enough to check that,
for every $\delta>0$,
$$\lim_{\varepsilon\to 0} \varepsilon^{-4}\,
\N_0\Big((1-e^{-\sigma})\int_0^\sigma ds\,
|F(W^{[s]})-F(W^{[s_*]})|\,
\ind{\{\wh W_s<\un W+\varepsilon,\;\zeta_s>\delta\}}\Big)=0.$$
However,
\ba
&& \varepsilon^{-4}\,
\N_0\Big((1-e^{-\sigma})\int_0^\sigma ds\,
|F(W^{[s]})-F(W^{[s_*]})|\,
\ind{\{\wh W_s<\un W+\varepsilon,\;\zeta_s>\delta\}}\Big)\\
&&\ \leq  \varepsilon^{-4}\,
\N_0\Big((1-e^{-\sigma})
\sup_{\{s\in[0,\sigma]:\wh W_s<\un W+\varepsilon\}}
(|F(W^{[s]})-F(W^{[s_*]})|)\;
\int_0^\sigma ds\,
\ind{\{\wh W_s<\un W+\varepsilon,\;\zeta_s>\delta\}}\Big)\\
&&\ \leq  C_\delta\,
\N_0\Big((1-e^{-\sigma})^2
\sup_{\{s\in[0,\sigma]:\wh W_s<\un W+\varepsilon\}}
(|F(W^{[s]})-F(W^{[s_*]})|^2)\Big)^{1/2}
\ea
by the Cauchy-Schwarz inequality and Lemma \ref{exccondlem2}. The last 
quantity tends to $0$ by (\ref{exccond1}) and dominated convergence. \cq

From (\ref{reroot1}) and Corollary \ref{exccondcoro}, the convergence
of Theorem \ref{exc-cond} reduces to checking that
\be
\label{exccondkey}
\lim_{\varepsilon\to 0}
\varepsilon^{-4}\;\N_0\Big(\varphi(\sigma)\,
F(W^{[s_*]})
\int_0^\sigma ds\,
\ind{\{\wh W_s<\un W+\varepsilon\}}\Big)
={2\over 21}\;\N_0(\varphi(\sigma)\,F(W^{[s_*]})).
\ee
The proof of (\ref{exccondkey}) will require two more lemmas. Before
stating
the first one, we need to introduce some notation. For $x\geq 0$, we
suppose that we are given a Poisson point measure
$\n=\sum_{i\in I}\delta_{\omega^i}$ with intensity $x\N_0$,
under the probability measure $\P_{(x)}$. To simplify notation,
we write $W^i_s=W_s(\omega^i)$. We then set
$$\un{\un W}=\inf_{i\in I}\Big(\inf_{s\geq 0} \wh W^i_s\Big)=\inf _{i\in
I}
\un W^i.$$

\begin{lemma}
\label{expliconst}
For every $x>0$ and $\varepsilon>0$,
$$\E_{(x)}\Big[\sum_{i\in I} \int_0^{\sigma(\omega^i)} ds\,
\ind{\{\wh W^i_s<\un{\un W}+\varepsilon\}}\Big]
=\varepsilon^4\,g({\varepsilon
\over \sqrt{x}}),$$
where the function $g:\R_+\la \R_+$ is continuous and nonincreasing,
and $g(0)=2/21$. 
\end{lemma}

\proof We first recall a well-known fact about Palm distributions of Poisson point measures.
If $\m$ is a Poisson point measure on a locally compact space, and if the intensity measure $m$
of $\m$ is a Radon measure, then, for every nonnegative measurable functional 
$\Phi$,
$$E\Big[\int \m(de)\,\Phi(e,\m)\Big]=
\int m(de)\,E[\Phi(e,\m+\delta_e)].$$
See e.g. Sections 10 and 11 in \cite{Kal}. We apply this to the point measure 
$\n$ and to the function
$$\Phi(\omega,\n)=\int_0^{\sigma(\omega)} ds \,\ind{\{\wh W_s(\omega)<\un{\un W}+\varepsilon\}}.$$ 
Note that
$\Omega$ is not locally compact, but as a Polish space it is homeomorphic to a Borel subset of a compact
metric space, so that the application of the preceding formula is easy to justify in our setting. We get
\ba
\E_{(x)}\Big[\sum_{i\in I} \int_0^{\sigma(\omega^i)} ds\,
\ind{\{\wh W^i_s<\un{\un W}+\varepsilon\}}\Big]
&=&x\,\N_0\Big(\int_0^\sigma ds\,\ind{\{\wh W_s<\un W+\varepsilon\}}
\P_{(x)}[\un{\un W}>a-\varepsilon]_{a=\wh W_s}\Big)\\
&=&x\,\N_0\Big(\int_0^\sigma ds\,\ind{\{\wh W_s<\un W+\varepsilon\}}
\,\exp\Big(-x\N_0(\un W\leq a-\varepsilon)_{a=\wh W_s}\Big)\Big)\\
&=&x\,\N_0\Big(\int_0^\sigma ds\,\ind{\{\wh W_s<\un W+\varepsilon\}}
\exp\Big(-{3x\over 2(\hat W_s-\varepsilon)^2}\Big)\Big)
\ea
by Lemma \ref{hitting}(i). In a way analogous to the proof of Lemma
\ref{exccondlem1} above, this quantity is equal to
\ba 
&&x\int_0^\infty
da\,E_0\Big[\ind{\{\xi[0,a]\subset]\xi_a-\varepsilon,\infty[\}}
\;\exp\Big(-6\int_0^a {dt\over (\xi_t-\xi_a+\varepsilon)^2}
-{3x\over 2(\xi_a-\varepsilon)^2}\Big)\Big]\\
&&\qquad=
x\int_0^\infty
da\,E_\varepsilon\Big[\ind{\{\xi[0,a]\subset]0,\infty[\}}
\;\exp\Big(-6\int_0^a {dt\over \xi_t^2}
-{3x\over 2\xi_a^2}\Big)\Big]\\
&&\qquad =
x\,\varepsilon^4\int_0^\infty
da\,E^{(9)}_\varepsilon\Big[R_a^{-4}\exp-{3x\over 2R_a^2}\Big]\\
&&\qquad =\varepsilon^4\int_0^\infty
db\,E^{(9)}_{\varepsilon/\sqrt{x}}\Big[R_b^{-4}\exp-{3\over 2R_b^2}\Big],
\ea
which gives the formula of the lemma, with
$$g(u)=E^{(9)}_u\Big[\int_0^\infty db\,R_b^{-4}\exp-{3\over
2R_b^2}\Big].$$
The fact that $g$ is nondecreasing follows from the strong Markov
property of $R$. The continuity of $g$ is easy from
a similar argument. Finally, the value of $g(0)$ is obtained from
the explicit formula for the Green function of
nine-dimensional Brownian motion. \cq

Recall our notation $\e^{]a,\infty[}$ for the $\sigma$-field generated
by the Brownian snake paths before their first exit from $]a,\infty[$,
and $L^a$ for the total mass of the exit measure $\z^{]a,\infty[}$. 

\begin{lemma}
\label{spemar}
Let $a<0$. For every bounded $\e^{]a,\infty[}$-measurable function 
$\Phi$ on $\Omega$,
$$\lim_{\varepsilon\to 0}
\N_0\Big(\ind{\{\un W<a\}}\,\Phi\,\int_0^\sigma ds\,\ind{\{\wh W_s<\un
W+\varepsilon\}}\Big)={2\over 21}\;\N_0(\ind{\{\un W<a\}}\,\Phi).$$
\end{lemma}

\proof For every $\w\in\W_0$, set $\tau_a(\w)=\inf\{t\geq 0:\w(t)\leq
a\}$. We first show that
\be 
\label{spemar1}
\lim_{\varepsilon\to 0}
\varepsilon^{-4}\;\N_0\Big(\ind{\{\un W<a\}}\,\int_0^\sigma ds\,\ind{\{\wh
W_s<\un W+\varepsilon,\;\tau_a(W_s)\geq\zeta_s\}}\Big)=0.
\ee
Clearly, it is enough to prove that
\be 
\label{spemar2}
\lim_{\varepsilon\to 0}
\varepsilon^{-4}\;\N_0\Big(\ind{\{\un W<a\}}\,\int_0^\sigma ds\,\ind{\{\wh
W_s<\un W+\varepsilon,\;\wh W_s\geq a\}}\Big)=0.
\ee
If $\delta>0$ is fixed, we have first
\ba
\varepsilon^{-4}\;\N_0\Big(\ind{\{\un W<a\}}\,\int_0^\sigma ds\,\ind{\{\wh
W_s<\un W+\varepsilon,\;\wh W_s\geq a\}}\ind{\{\zeta_s>\delta\}}\Big)
&\leq&
\varepsilon^{-4}\;\N_0\Big(\ind{\{a-\varepsilon<\un W<a\}}\,\int_0^\sigma
ds\,\ind{\{\wh W_s<\un W+\varepsilon,\;\zeta_s>\delta\}}\Big)\\
&\leq&C_\delta\,\N_0(a-\varepsilon<\un W<a),
\ea
by the Cauchy-Schwarz inequality and Lemma \ref{exccondlem2}. Obviously
the last quantity tends to $0$ as $\varepsilon\to 0$. Then,
\ba
&&\varepsilon^{-4}\;\N_0\Big(\ind{\{\un W<a\}}\,
\int_0^\sigma ds\,\ind{\{\wh
W_s<\un W+\varepsilon,\;\wh W_s\geq a\}}\ind{\{\zeta_s\leq\delta\}}\Big)\\
&&\qquad\leq 
\varepsilon^{-4}\;\N_0\Big(
\int_0^\sigma ds\,\ind{\{\wh
W_s<(\un
W+\varepsilon)\wedge(a+\varepsilon)\}}
\,\ind{\{\zeta_s\leq\delta\}}\Big)\\
&&\qquad= 
\varepsilon^{-4}\int_0^\delta dt\,E_0\Big[
\ind{\{\xi_t<a+\varepsilon,\;\xi[0,t]\subset]\xi_t-\varepsilon,\infty[\}}
\,\exp\Big(-6\int_0^t {dr\over(\xi_r-\xi_t+\varepsilon)^2}\Big)\Big]\\
&&\qquad=
\varepsilon^{-4}\int_0^\delta dt\,E_\varepsilon\Big[
\ind{\{\xi_t>-a,\;\xi[0,t]\subset]0,\infty[\}}
\,\exp\Big(-6\int_0^t {dr\over\xi_r^2}\Big)\Big]\\
&&\qquad=\int_0^\delta dt\,E^{(9)}_\varepsilon[\ind{\{R_t>-a\}}\,R_t^{-4}]
\\
&&\qquad\leq \delta\,a^{-4}\;.
\ea
The last quantity can be made arbitrarily small by choosing $\delta$
small, independently of $\varepsilon$. This completes the 
proof of (\ref{spemar2}) and (\ref{spemar1}). 

It remains to study
\ba
&&\varepsilon^{-4}\;\N_0\Big(\ind{\{\un W<a\}}\,\Phi\,\int_0^\sigma
ds\,\ind{\{\wh W_s<\un W+\varepsilon,\;\tau_a(W_s)<\zeta_s\}}\Big)\\
&&\qquad = \varepsilon^{-4}\;\N_0\Big(\ind{\{\un W<a\}}\,\Phi\,
\E_{(L^a)}
\Big[\sum_{i\in I} \int_0^{\sigma(\omega^i)} ds\,
\ind{\{\wh W^i_s<\un{\un W}+\varepsilon\}}\Big]\Big)\\
&&\qquad=\N_0\Big(\ind{\{\un W<a\}}\,\Phi\,
g\Big({\varepsilon\over \sqrt{L^a}}\Big)\Big).
\ea
In the first equality we used the special Markov property (Theorem
\ref{SMP}), and in the second one Lemma \ref{expliconst}. The desired
result now follows from dominated convergence. \cq

\noindent{\bf Proof of Theorem \ref{exc-cond}:} \ We already noticed that
it
is enough to establish (\ref{exccondkey}). We first observe that
\be
\label{exckey1}
\lim_{b\da 0}
\sup_{0<\varepsilon<1} \left(\varepsilon^{-4}\,
\N_0\Big(\ind{\{\un W\geq -b\}}\,\varphi(\sigma)\int_0^\sigma ds\,
\ind{\{\wh W_s<\un W+\varepsilon\}}\Big)\right)=0.
\ee
In fact, for every $\delta>0$, Lemma \ref{exccondlem2}
gives the bound
$$\varepsilon^{-4}\,
\N_0\Big(\ind{\{\un W\geq -b\}}\,\varphi(\sigma)\int_0^\sigma ds\,
\ind{\{\wh W_s<\un W+\varepsilon,\;\zeta_s>\delta\}}\Big)
\leq C_\delta\,\N_0\Big(1_{\{\un W\geq -b\}}\,\varphi(\sigma)^2
\Big)^{1/2},$$
and the right-hand side tends to $0$ as $b\da 0$ by dominated convergence.
On the other hand, 
$$\sup_{0<\varepsilon<1} \left(\varepsilon^{-4}\,
\N_0\Big(\varphi(\sigma)\int_0^\sigma ds\,
\ind{\{\wh W_s<\un W+\varepsilon,\;\zeta_s\leq \delta\}}\Big)\right)\;
\build{\la}_{\delta\da 0}^{}\; 0$$
by Lemma \ref{exccondlem1}. This completes the proof
of (\ref{exckey1}).

For every $a<0$, set 
$$T_a=\inf\{s\geq 0:\wh W_s=a\},$$
and write 
$\wt W^a:=\wt W^{]a,\infty[}$ for the ``Brownian snake truncated
below level $a$'' (cf subsection 2.4). By definition, $\wt W^a$
is $\e^{]a,\infty[}$-measurable. Furthermore, we have also
$$T_a=\inf\{s\geq 0:\wh{\wt W}_s=a\},\quad \N_0\hbox{ a.e.}$$
For every $s\in[0,\sigma(\wt W^a)]$, we can define the re-rooted snake
$\wt W^{a,[s]}$ from $\wt W^a$,
in the same way as $W^{[s]}$ was defined from
$W$ in Section 2. Note that the process 
$\wt W^{a,[T_a]}$, whose definition makes sense
on the $\e^{]a,\infty[}$-measurable set $\{T_a<\infty\}$, is 
$\e^{]a,\infty[}$-measurable. Hence, Lemma \ref{spemar} gives
\begin{eqnarray}
\label{exckey2}
&&\lim_{\varepsilon\to 0}
\varepsilon^{-4}\;\N_0\Big(\ind{\{\un W<a\}}\,\varphi(\sigma(\wt W^a))\,
F(\wt W^{a,[T_a]})
\int_0^\sigma ds\,
\ind{\{\wh W_s<\un W+\varepsilon\}}\Big)\nonumber\\
&&\qquad = {2\over 21}\;\N_0\Big(\ind{\{\un W<a\}}\,\varphi(\sigma(\wt W^a))\,
F(\wt W^{a,[T_a]})\Big).
\end{eqnarray}

We fix $a<0$ of the form $a=k_0\,2^{-n_0}$, where $k_0\in\Z$ and $n_0\in
\N$. For every $x<0$ and $n\in \N$, we denote by $\{x\}_n$ the smallest
number
of the form $k\,2^{-n}$, with $k\in[-2^{2n},0]\cap \Z$, which is strictly
greater than $x$. As a consequence of Lemma \ref{spemar} (applied 
with $a$ replaced by $i2^{-n}$ and with a suitable choice of $\Phi$),
we get, for every integer $n\geq n_0$ and every $i\in \Z$ such that
$-2^{n}\leq i2^{-n}<a$,
\begin{eqnarray}
\label{exckey3}
&&\lim_{\varepsilon\da 0} 
\varepsilon^{-4}\, \N_0\Big(\ind{\{\un W<i2^{-n}\}}
\Big(\varphi(\sigma(\wt W^{i
2^{-n}}))\,F(
\wt W^{i2^{-n},[T_{i2^{-n}}]})\nonumber\\
&&\hspace{4cm}-\varphi(\sigma(\wt W^{(i+1)
2^{-n}}))\,F(
\wt W^{(i+1)2^{-n},[T_{(i+1)2^{-n}}]})\Big) \int_0^\sigma ds\,\ind{\{
\wh W_s<\un W+\varepsilon\}}\Big)\nonumber\\
&&={2\over 21}\; \N_0\Big(\ind{\{\un W<i2^{-n}\}}
\Big(\varphi(\sigma(\wt W^{i
2^{-n}}))\,F(
\wt W^{i2^{-n},[T_{i2^{-n}}]})\nonumber\\
&&\hspace{4cm}-\varphi(\sigma(\wt W^{(i+1)
2^{-n}}))\,F(
\wt W^{(i+1)2^{-n},[T_{(i+1)2^{-n}}]})\Big)\Big).
\end{eqnarray}
We sum (\ref{exckey2}) and the convergences (\ref{exckey3}) for all
choices of $i\in\Z$ with $-2^n\leq i\,2^{-n}<a$. It follows that
\begin{eqnarray}
\label{exckey4}
&&\lim_{\varepsilon\da 0} 
\varepsilon^{-4}\, \N_0\Big(\ind{\{\un W<a\}}
\varphi(\sigma(\wt W^{\{\un W\}_n}
))\,F(
\wt W^{\{\un W\}_n,[T_{\{\un W\}_n}]})\int_0^\sigma ds\,\ind{\{
\wh W_s<\un W+\varepsilon\}}\Big)\nonumber\\
&&\qquad={2\over 21}\;\N_0\Big(\ind{\{\un W<a\}}
\varphi(\sigma(\wt W^{\{\un W\}_n}
))\,F(
\wt W^{\{\un W\}_n,[T_{\{\un W\}_n}]})\Big).
\end{eqnarray}
Note that $\sigma(\wt W^{\{\un W\}_n})\leq \sigma$, and
$$\sigma(\wt W^{\{\un W\}_n})\build{\la}_{n\to \infty}^{} \sigma\;,
\quad \N_0\hbox{ a.e.}$$
Moreover, $T_{\{\un W\}_n}\la T_{\un W}=s_*$, and from the construction
of re-rooted snakes it follows that
$$
\wt W^{\{\un W\}_n,[T_{\{\un W\}_n}]}
\build{\la}_{n\to \infty}^{} W^{[s_*]},$$
$\N_0$ a.e., in the sense of uniform convergence. By dominated
convergence,
we get that the right-hand side of (\ref{exckey4}) is close to
$${2\over 21}\;\N_0\Big(\ind{\{\un W<a\}}\,\varphi(\sigma)\,
F(W^{[s_*]})\Big)$$
when $n\to \infty$.

Using (\ref{exckey1}) and (\ref{exckey4}), we see that the 
proof of 
(\ref{exccondkey}) will be complete if we can verify that
$$\sup_{0<\varepsilon<1}
\Big(\varepsilon^{-4}\,\N_0\Big(\ind{\{\un W<a\}}\,
|\varphi(\sigma(\wt W^{\{\un W\}_n}
))\,F(
\wt W^{\{\un W\}_n,[T_{\{\un W\}_n}]})-\varphi(\sigma)\,
F(W^{[s_*]})|
\int_0^\sigma ds\,\ind{\{
\wh W_s<\un W+\varepsilon\}}\Big)\Big)
$$
tends to $0$ as $n\to \infty$. This is easy by decomposing 
the set $\{\wh W_s<\un W+\varepsilon\}$ as
$$\{\wh W_s<\un W+\varepsilon,\;\zeta_s\leq \delta\}
\cup \{\wh W_s<\un W+\varepsilon,\;\zeta_s> \delta\}$$
and using Lemma \ref{exccondlem1} for the first term and 
Lemma \ref{exccondlem2}, together with the Cauchy-Schwarz inequality,
for the second one, as we did previously. This completes the
proof of (\ref{exccondkey}) and of Theorem \ref{exc-cond}. \cq

\noindent{\bf Proof of Theorems \ref{main-cond} and \ref{verwaat}:} \ 
Both Theorems \ref{main-cond} and \ref{verwaat} follow from the
convergence 
\be 
\label{excnorm}
\lim_{\varepsilon\to 0}
\varepsilon^{-4}\;\N^{(1)}_0\Big(F(W)\,
\ind{\{\un W>-\varepsilon\}}\Big)={2\over 21}\;
\N^{(1)}_0(F(W^{[s_*]})),
\ee
which holds for every bounded continuous function $F$ on
$\Omega=C(\R_+,\W)$ (take $F=1$ to recover the
first assertion of Theorem \ref{main-cond}). We will now derive
(\ref{excnorm}) from Theorem \ref{exc-cond}. 

For every $\lambda>0$, let us introduce the scaling operator
$\theta_\lambda$ defined on $\Omega$ by
$$\begin{array}{l}
\zeta_s\circ \theta_\lambda=\lambda^{1/2}\,\zeta_{s/\lambda}\;,\\
\noalign{\smallskip}
W_s\circ
\theta_\lambda(t)=\lambda^{1/4}\,W_{s/\lambda}(\lambda^{-1/2}t)\;.
\end{array}
$$
Note that, for every $r>0$, the image of $\N^{(r)}_0$
under $\theta_{1/r}$ is $\N^{(1)}_0$.

Let $\delta\in]0,1[$. It follows from Theorem \ref{exc-cond}
that the law of the pair $(\sigma,W)$ under
$$\mu_{\varepsilon,\delta} :=\varepsilon^{-4}
\N_0(\cdot \cap \{\un W>-\varepsilon,\;1-\delta<\sigma<1\})$$
converges weakly as $\varepsilon\to 0$ towards the law of 
$(\sigma,W^{[s_*]})$ under the measure $\mu_\delta$
having density $2/(21\sigma)$ with respect to 
$\N_0(\cdot\cap\{1-\delta<\sigma<1\})$. 

Since the mapping $(r,\omega)\to \theta_r\omega$ is continuous,
it follows that the law of $W\circ \theta_{1/\sigma}$
under $\mu_{\varepsilon,\delta}$ converges as $\varepsilon\to 0$
towards the law of $W^{[s_*]}\circ \theta_{1/\sigma}$ 
under $\mu_\delta$. Thus,
$$\lim_{\varepsilon\to 0}\mu_{\varepsilon,\delta}(F(W\circ \theta_{1/\sigma}))
=\mu_\delta(F(W^{[s_*]}\circ\theta_{1/\sigma})),$$
or equivalently
$$\lim_{\varepsilon\to 0} \varepsilon^{-4}
\N_0\Big(\ind{\{\un W>-\varepsilon,\;1-\delta<\sigma<1\}}\,
F(W\circ \theta_{1/\sigma})\Big)
=\N_0\Big({2\over 21\sigma}\,\ind{\{1-\delta<\sigma<1\}}\,
F(W^{[s_*]}\circ\theta_{1/\sigma})\Big).$$
Since the density of $\sigma$ under $\N_0$
is $(8\pi)^{-1/2}\,s^{-3/2}$, this can be rewritten as
$$\lim_{\varepsilon\to 0}
\varepsilon^{-4}\int_{1-\delta}^1 {dr\over2\sqrt{2\pi r^3}}
\;\N^{(r)}_0\Big(F(W\circ \theta_{1/r})\,\ind{\{\un W >-\varepsilon\}}\Big)=
{2\over 21}\int_{1-\delta}^1 {dr\over2\sqrt{2\pi r^5}}
\;\N^{(r)}_0(F(W^{[s_*]}\circ \theta_{1/r})).$$
Now observe that $\un W=r^{1/4}
\,\un W\circ \theta_{1/r}$, and recall that the image of $\N^{(r)}_0$
under $\theta_{1/r}$ is $\N^{(1)}_0$ to get
$$
\lim_{\varepsilon\to 0}
\varepsilon^{-4}\int_{1-\delta}^1 {dr\over2\sqrt{2\pi r^3}}
\;\N^{(1)}_0\Big(F(W)\,\ind{\{\un W >-r^{-1/4}\varepsilon\}}\Big)
= {2\over 21}\Big(\int_{1-\delta}^1 {dr\over2\sqrt{2\pi r^5}}\Big)
\;\N^{(1)}_0(F(W^{[s_*]})).
$$
Without loss of generality we can assume that $F\geq 0$. Since
$r^{-1/4}\varepsilon\geq \varepsilon$ for $1-\delta\leq r\leq 1$, taking
$\delta$ small
in the preceding convergence leads to
$$\limsup_{\varepsilon\to 0} \varepsilon^{-4}\,
\N^{(1)}_0\Big(F(W)\,\ind{\{\un W >-\varepsilon\}}\Big)
\leq {2\over 21}\,\N^{(1)}_0(F(W^{[s_*]})).
$$
By arguing with the constraint $1<\sigma<1+\delta$ instead
of $1-\delta<\sigma<1$, we get the analogous lower bound for the
liminf behavior. This completes the proof of (\ref{excnorm}). \cq

We conclude this section with another approximation of the conditioned
measure $\ov\N^{(1)}_0$, which is similar to Theorem 
\ref{main-cond} but much easier to obtain.

\begin{proposition}
\label{indepunif}
Let $U_1,U_2,\ldots$ be a sequence of i.i.d. uniform $[0,1]$ 
random variables defined under an auxiliary probability measure $Q$.
Then, for any bounded continuous function $F$ on $\Omega$,
$$\lim_{p\to\infty} \N^{(1)}_0\otimes Q(F(W)\mid
\wh W_{U_1}>0,\wh W_{U_2}>0,\ldots,\wh W_{U_p}>0)=\ov \N^{(1)}_0(F).$$
\end{proposition}

\proof From the re-rooting theorem (Theorem \ref{reroot}) and the remark following this statement,
\ba
&&\N^{(1)}_0\otimes Q\Big(F(W)\,\ind{\{
\wh W_{U_1}>0,\ldots,\wh W_{U_p}>0\}}\Big)\\
&&\qquad=\N^{(1)}_0\otimes Q\Big(\int_0^1 ds\,F(W^{[s]})\,
\ind{\{
\wh W^{[s]}_{U_1}>0,\ldots,\wh W^{[s]}_{U_p}>0\}}\Big)\\
&&\qquad=\N^{(1)}_0\otimes Q\Big(\int_0^1 ds\,F(W^{[s]})\,
\ind{\{
\wh W_{s\oplus U_1}>\wh W_s,\ldots,\wh W_{s\oplus U_p}>\wh W_s\}}\Big)\\
&&\qquad=\N^{(1)}_0\Big(
\int_0^1du_0\int_0^1 du_1 \ldots\int_0^1 du_p\,
F(W^{[u_0]})\,
\ind{\{
\wh W_{u_1}>\wh W_{u_0},\ldots,\wh W_{u_p}>\wh W_{u_0}\}}\Big)\\
&&\qquad=
{1\over p+1}\sum_{i=0}^{p} \N^{(1)}_0\Big(
\int_0^1du_0\int_0^1 du_1 \ldots\int_0^1 du_p\,
F(W^{[u_i]})\,
\ind{\{
\wh W_{u_j}>\wh W_{u_i},\;\forall j\not= i\}}\Big)\\
&&\qquad=
{1\over p+1}\;\N^{(1)}_0\Big(
\int_0^1du_0\int_0^1 du_1 \ldots\int_0^1 du_p\,
F(W^{[u^{(p)}_{\rm min}]})\Big)
\ea
where $u^{(p)}_{\rm min}:=u_{j^{(p)}_{\rm min}}$, if $j^{(p)}_{\rm min}$
is the (a.e. unique) index $i$ such that $\wh W_{u_i}=\inf\{\wh W_{u_j},\,
1\leq j\leq p\}$. Taking $F=1$, we have
$$\N^{(1)}_0\otimes Q(
\wh W_{U_1}>0,\ldots,\wh W_{U_p}>0)={1\over p+1}$$
and, on the other hand, dominated convergence shows that
$$\lim_{p\to\infty} 
\N^{(1)}_0\Big(
\int_0^1du_0\int_0^1 du_1 \ldots\int_0^1 du_p\,
F(W^{[u^{(p)}_{\rm min}]})\Big)
=\N^{(1)}_0(F(W^{[s_*]}))=\ov\N^{(1)}_0(F).$$
This completes the proof of Proposition \ref{indepunif}. \cq

\section{Other conditionings}

Motivated by Theorem \ref{exc-cond}, we define a $\sigma$-finite measure 
$\ov \N_0$ on $\Omega$ by setting
$$\ov \N_0(F)=\N_0({1\over \sigma}\,F(W^{[s_*]})).$$
Theorem \ref{exc-cond} shows that, up to the multiplicative constant
$2/21$, $\ov \N_0$ is the limit in an appropriate sense of the
measures $\varepsilon^{-4}\,\N_0(\cdot \cap \{\un W>-\varepsilon\})$
as $\varepsilon\to 0$. We have also
$$\ov \N_0(F)=\int_0^\infty {dr\over 2\sqrt{2\pi r^5}}\,
\N_0^{(r)}(F(W^{[s_*]}))
=\int_0^\infty {dr\over 2\sqrt{2\pi r^5}}\,\ov \N_0^{(r)}(F),$$
where $\ov \N_0^{(r)}$ can be defined equivalently as the
law of $W^{[s_*]}$ under $\N^{(r)}_0$, or
as the image of $\ov\N^{(1)}_0$ under the scaling operator $\theta_r$.

We will now describe a different approach to $\ov \N_0$, which involves
conditioning the Brownian snake excursion on its height $H=\sup_{s\geq 0}
\zeta_s$, rather than on its length as in Theorem \ref{main-cond}. 
This will give more insight in the behavior of the Brownian snake
under $\ov\N_0$. Eventually, this will lead to a construction of
a Brownian snake excursion with infinite length conditioned to
stay on the positive side. We rely on some ideas from \cite{AW}.

For every $h>0$, we set $\N^h_0=\N_0(\cdot\mid H=h)$. Then,
$$\N_0=\int_0^\infty {dh\over 2h^2}\;\N^h_0.$$
From Theorem 1 in \cite{AW} we know that there exists a constant
$c_0>0$ such that
\be 
\label{AbW}
\lim_{\varepsilon\to 0} \varepsilon^{-4}\,\N^1_0(\un W>-\varepsilon)=c_0.
\ee
A simple scaling argument then implies that, for every $h>0$,
\be 
\label{AbW2}
\lim_{\varepsilon\to 0} \varepsilon^{-4}\,\N^h_0(\un W>-\varepsilon)=
{c_0\over h^2}.
\ee

\begin{theorem}
\label{heightcond}
For every $h>0$, there exists a probability measure $\ov\N^h_0$ on $\Omega$
such that 
$$\lim_{\varepsilon\to 0} \N^h_0(\cdot \mid \un
W>-\varepsilon)=\ov\N^h_0$$
in the sense of weak convergence on the space of probability measures on
$\Omega$. Moreover,
$$\ov\N_0={21 c_0\over 4} \int_0^\infty {dh\over h^4}\,\ov\N^h_0.$$
\end{theorem}

\noindent{\bf Remark.} Our proof of the first part of Theorem
\ref{heightcond} does not use Section 3. This proof thus gives
another approach to the conditioned measure $\ov\N_0$, which does
not depend on the re-rooting method that played a crucial role
in Section 3.

Before proving Theorem \ref{heightcond}, we will establish an
important preliminary result. We first introduce some notation.
Following \cite{AW}, we set for every $\varepsilon>0$,
$$f(\varepsilon)=\N^1_0(\un W>-\varepsilon)$$
and, for every $x>0$,
$$G(x)=4\int_0^x u(1-f(u))\,du.$$
The function $G$ is obviously nondecreasing. It is also bounded
since
$$
G(\infty)=4\int_0^\infty u\,\N^1_0(\un W\leq -u)\,du
=2\int_0^\infty r^{-2}\,\N^r_0(\un W\leq -1)\,dr=4\,\N_0(\un W\leq -1) =
6$$
by a scaling argument and Lemma \ref{hitting}(i).

By well-known properties of Brownian excursions, there exists 
$\N^h_0$ a.s. a unique time $\alpha\in]0,\sigma[$ such that
$\zeta_\alpha=h$. The next proposition discusses the law of
$W_\alpha$ under $\N^h_0(\cdot\mid \un W>-\varepsilon)$. 

\begin{proposition}
\label{heighttech}
Let $\Phi$ be a bounded continuous function on $\W$. Then,
$$\lim_{\varepsilon\da 0} \varepsilon^{-4}\,\N^h_0(\Phi(W_\alpha)
\,\ind{\{\un W>-\varepsilon\}})
=E^{(9)}_0\Big[\Phi(R_t,0\leq t\leq h)\,R_h^{-4}\,\exp
\Big(\int_0^h{dt\over R_t^2}\,G\Big({R_t\over
\sqrt{h-t}}\Big)\Big)\Big].$$
\end{proposition}

\noindent{\bf Remarks.} (i) From the bound $G(x)\leq 6\wedge (2x^2)$, it
is immediate to verify that
$$\int_0^h{dt\over R_t^2}\,G\Big({R_t\over
\sqrt{h-t}}\Big)<\infty\;,\qquad P^{(9)}_0\hbox{ a.s.}$$
(ii) By taking $\Phi=1$, we see that the constant $c_0$ in 
(\ref{AbW}) is given by
$$c_0=E^{(9)}_0\Big[R_1^{-4}\,\exp
\Big(\int_0^1{dt\over R_t^2}\,G\Big({R_t\over
\sqrt{1-t}}\Big)\Big)\Big],$$
as it was already observed in \cite{AW}. The fact that the quantity in
the right-hand side is finite follows from the proof below.

\proof Our main tool is Williams' decomposition of the Brownian excursion
at its maximum (see e.g. Theorem XII.4.5 in \cite{RY}). For every $s\geq
0$, we set
$$\rho_s=\zeta_{s\wedge \alpha}\ ,\quad \rho'_s=\zeta_{(\sigma-s)\vee
\alpha}.$$
Under the probability measure $\N^h_0$, the processes $(\rho_s)_{s\geq 0}$
and $(\rho'_s)_{s\geq 0}$ are two independent three-dimensional Bessel
processes started at $0$ and stopped at their first hitting time of $h$. 

We also need to introduce the excursions of $\rho$ and $\rho'$ above their
future infimum. Set
$$\un \rho_s=\inf_{r\geq s} \rho_r$$
and let $(a_j,b_j),\;j\in J$ be the connected components of the
open set $\{s\geq 0:\rho_s>\un \rho_s\}$. For every $j\in J$, define
\ba
&&\zeta^j_s=\rho_{(a_j+s)\wedge b_j}-\rho_{a_j}\;,\qquad s\geq 0\\
&&h_j=\rho_{a_j}\;.
\ea
Then, by excursion theory,
$$\sum_{j\in J} \delta_{(h_j,\zeta^j)}(dr\,de)$$
is a Poisson point measure on $\R_+\times C(\R_+,\R_+)$ with 
intensity
$$2\;\ind{[0,h]}(r)\;\ind{[0,h-r]}(H(e))\,dr\,n(de)$$
where $H(e)=\sup_{s\geq 0}e(s)$ as previously. The same result obviously
holds for the analogous point measure
$$\sum_{j\in J'} \delta_{(h'_j,\zeta'^j)}(dr\,de)$$
obtained by replacing $\rho$ with $\rho'$.

We can combine the preceding assertions with the spatial
displacements of the Brownian snake, in a way very similar to the proof
of Lemma V.5 in \cite{Zu}. For every
$j\in J$, we set
$$W^j_s(t)=W_{(a_j+s)\wedge b_j}(h_j+t)-\wh W_{a_j}\;,\quad 0\leq t\leq
\zeta^j_s,\quad s\geq 0.$$
Note that by the properties of the Brownian snake $\wh W_{a_j}=
\wh W_{b_j}=W_\alpha(h_j)$. Then,
$$\n:=\sum_{j\in J} \delta_{(h_j,W^j)}$$
is under $\N^h_0$ a Poisson point measure on $\R_+\times \Omega$, with intensity
\be
\label{Poisson1}
2\;\ind{[0,h]}(r)\;\ind{[0,h-r]}(H(\omega))\,dr\,\N_0(d\omega).
\ee
The same holds for the analogous point measure
$$\n':=\sum_{j\in J'} \delta_{(h'_j,W'^j)}.$$
Moreover $\n$ and $\n'$ are independent and the pair $(\n,\n')$
is independent of $W_\alpha$. All these assertions easily follow from
properties of the Brownian snake.

Now note that the range of the Brownian snake under $\N^h_0$
can be decomposed as
$$\{W_\alpha(t):0\leq t\leq h\}
\cup \Big(\bigcup_{j\in J} (W_\alpha(h_j)+\r(W^j))\Big)
\cup \Big(\bigcup_{j\in J'} (W_\alpha(h'_j)+\r(W'^j))\Big).$$
Using this observation and conditioning with respect to $W_\alpha$,
we get
\ba 
&&\N^h_0(\Phi(W_\alpha)
\,\ind{\{\un W>-\varepsilon\}})\\
&&\quad =\N^h_0\Big(\Phi(W_\alpha)
\,\ind{\{W_\alpha(t)>-\varepsilon,\;0\leq t\leq h\}}
\exp\Big(-4\int_0^h dt\,\N_{W_\alpha(t)}(H<h-t,\,\un W\leq -\varepsilon)
\Big)\Big)\\
&&\quad =E_0\Big[\Phi(\xi_t,\,0\leq t\leq h)
\,\ind{\{\xi[0,h]\subset]-\varepsilon,\infty[\}}
\exp\Big(-4\int_0^h dt\,\N_{\xi_t}(H<h-t,\,\un W\leq -\varepsilon)
\Big)\Big].
\ea
Then, for every $x>-\varepsilon$,
$$\N_{x}(H<h-t,\,\un W\leq -\varepsilon)
=\int_0^{h-t} {du\over 2u^2}\,\N^u_x(\un W\leq -\varepsilon)
=\int_0^{h-t} {du\over 2u^2}\,(1-f({x+\varepsilon\over \sqrt{u}}))$$
and we obtain
\begin{eqnarray}
\label{heighttech1}
&&\N^h_0(\Phi(W_\alpha)
\,\ind{\{\un W>-\varepsilon\}})\nonumber\\
&&\ =E_0\Big[\Phi(\xi_t,\,0\leq t\leq h)
\,\ind{\{\xi[0,h]\subset]-\varepsilon,\infty[\}}
\exp\Big(-2\int_0^h dt
\int_0^{h-t} {du\over u^2}\,(1-f({\xi_t+\varepsilon\over \sqrt{u}}))\Big)
\Big]\nonumber\\
&&\ =E_\varepsilon\Big[\Phi(\xi_t-\varepsilon,\,0\leq t\leq h)
\,\ind{\{\xi[0,h]\subset]0,\infty[\}}
\exp\Big(-2\int_0^h dt
\int_0^{h-t} {du\over u^2}\,(1-f({\xi_t\over \sqrt{u}}))\Big)
\Big].
\end{eqnarray}
For every $x>0$, the change of variable $v=x/\sqrt{u}$ gives
$$\int_0^{h-t} {du\over u^2}\,(1-f({x\over \sqrt{u}}))
=2x^{-2}\int^\infty_{x/\sqrt{h-t}} dv\,v(1-f(v))
=x^{-2}\Big(3-{1\over 2}\,G\Big({x\over \sqrt{h-t}}\Big)\Big).$$
By substituting this into (\ref{heighttech1}) and using
Proposition \ref{Bessel} once more, we get
\begin{eqnarray}
\label{heighttech2}
&&\N^h_0(\Phi(W_\alpha)
\,\ind{\{\un W>-\varepsilon\}})\nonumber\\
&&\ =E_\varepsilon\Big[\Phi(\xi_t-\varepsilon,\,0\leq t\leq h)
\,\ind{\{\xi[0,h]\subset]0,\infty[\}}
\exp\Big(-6\int_0^h {dt\over \xi_t^2} +\int_0^h {dt\over \xi_t^2}\,
G\Big({\xi_t\over \sqrt{h-t}}\Big)\Big)\Big]\nonumber\\
&&\ =\varepsilon^4\,E^{(9)}_\varepsilon\Big[\Phi(R_t-\varepsilon,0\leq
t\leq h)\,R_h^{-4}\,\exp\Big(\int_0^h {dt\over R_t^2}\,
G\Big({R_t\over \sqrt{h-t}}\Big)\Big)\Big].
\end{eqnarray}
In view of (\ref{heighttech2}), the proof of Proposition \ref{heighttech}
reduces to checking that
\begin{eqnarray}
\label{heighttech3}
&&\lim_{\varepsilon\da 0}\;
E^{(9)}_\varepsilon\Big[\Phi(R_t-\varepsilon,0\leq
t\leq h)\,R_h^{-4}\,\exp\Big(\int_0^h {dt\over R_t^2}\,
G\Big({R_t\over \sqrt{h-t}}\Big)\Big)\Big]\nonumber\\
&&\quad= E^{(9)}_0\Big[\Phi(R_t,0\leq
t\leq h)\,R_h^{-4}\,\exp\Big(\int_0^h {dt\over R_t^2}\,
G\Big({R_t\over \sqrt{h-t}}\Big)\Big)\Big]
\end{eqnarray}
This follows from a dominated convergence argument, which at the same time
will prove that the quantity in the right-hand side of 
(\ref{heighttech3}) is well-defined. Note that we may define
on a common probability space, a nine-dimensional Bessel process
$X^\varepsilon=(X^\varepsilon_t,t\geq 0)$ started at $\varepsilon$,
for every $\varepsilon\geq 0$, in such a way that the inequality 
$X^\varepsilon\geq X^0$ 
holds a.s. for every $\varepsilon>0$. Since
$$G\Big({X^\varepsilon_t\over\sqrt{h-t}}\Big)\leq
{4(X^\varepsilon_t)^2\over h}\ ,
\quad \forall t\in[0,h/2]$$
we first get
\begin{eqnarray}
\label{heighttech4}
(X^\varepsilon_h)^{-4}\,\exp\Big(\int_0^h {dt\over
(X^\varepsilon_t)^2}\, G\Big({X^\varepsilon_t\over \sqrt{h-t}}\Big)\Big)
&\leq& e^2\,(X^\varepsilon_h)^{-4}\,\exp\Big(\int_{h/2}^h {dt\over
(X^\varepsilon_t)^2}\, G\Big({X^\varepsilon_t\over
\sqrt{h-t}}\Big)\Big)\nonumber\\ &\leq&
e^2\,(X^0_h)^{-4}\,\exp\Big(6\int_{h/2}^h {dt\over (X^0_t)^2}\Big),
\end{eqnarray}
using the bounds $G\leq 6$ and $X^\varepsilon\geq X^0$. Then, an
application of It\^o's formula shows that
$$(X^0_t)^{-4}\,\exp\Big(6\int_{h/2}^t {dr\over
(X^0_r)^2}\Big)$$
is a local martingale on the time interval $[h/2,\infty[$, and so
$$E\Big[(X^0_h)^{-4}\,\exp\Big(6\int_{h/2}^h {dt\over
(X^0_t)^2}\Big)\Big]\leq E[(X_{h/2}^0)^{-4}]<\infty.$$
Together with (\ref{heighttech4}), this shows that the random variables appearing in the left-hand side of
(\ref{heighttech3}) are uniformly integrable. The convergence (\ref{heighttech3})
easily follows. \cq

\noindent{\bf Proof of Theorem \ref{heightcond}:} \ 
We first explain how the first part of Theorem \ref{heightcond}
can be deduced from Proposition \ref{heighttech}. Recall the notation
$(\n,\n')$ from the proof of this proposition. We first observe that we
can find a measurable functional $\Gamma$ such that
$$W=\Gamma(W_\alpha,\n,\n')\;,\quad \N^h_0\hbox{ a.s.}$$
Let us make this functional more explicit. We have first 
$$\alpha=\sum_{j\in J} \sigma(W^j).$$
For every $\ell\in[0,h]$, we set
$$\tau_\ell=\sum_{j\in J} \ind{\{h_j\leq \ell\}}\,\sigma(W^j).$$
Then, if $s\in[0,\alpha]$, there is a unique $\ell$ such that
$\tau_{\ell-}\leq s\leq \tau_\ell$, and:
\begin{description}
\item{$\bullet$} Either there is a (unique) $j\in J$ such that $\ell=h_j$,
and
\ba
&&\zeta_s=\ell+\zeta^j_{s-\tau_{\ell-}}\;,\\
&&W_s(t)=\left\{\begin{array}{ll}
W_\alpha(t)&{\rm if}\ t\leq \ell\;,\\
W_\alpha(\ell)+W^j_{s-\tau_{\ell-}}(t-\ell)\qquad&{\rm if}\ \ell<t\leq
\zeta_s\;;
\end{array}
\right.
\ea
\item{$\bullet$} Or there is no such $j$, and 
\ba
&&\zeta_s=\ell\;,\\
&&W_s(t)=W_\alpha(t)\;,\qquad t\leq \ell\;.
\ea
\end{description}
The previous formulas identify $(W_s,0\leq s\leq \alpha)$ as a 
measurable function of the pair $(W_\alpha,\n)$, and in a similar way
we can recover $(W_{\sigma-s},0\leq s\leq \sigma-\alpha)$ as 
the same measurable function of $(W_\alpha,\n')$. 

To simplify notation, write $\N^{h,(\varepsilon)}$ for the conditional
probability $\N^h_0(\cdot\mid \un W>-\varepsilon)$. From elementary
properties of Poisson measures, we get that under the probability
measure $\N^{h,(\varepsilon)}_0$ and conditionally given $W_\alpha$,
the point measures $\n$ and $\n'$ are
independent and Poisson with intensity
$$\mu^h_\varepsilon(W_\alpha;dr\,d\omega):=2\;\ind{[0,h]}(r)\;\ind{[0,h-r]}(H(\omega))
\,\ind{\{\r(\omega)\subset]-\varepsilon-W_\alpha(r),\infty[\}}
\,dr\,\N_0(d\omega).$$
As a consequence of Proposition \ref{heighttech}, the law of $W_\alpha$
under
$\N^{h,(\varepsilon)}_0$
converges as
$\varepsilon\to 0$ to the law of the process $Y^h=(Y^h_t,0\leq t\leq h)$
such that
\be
\label{heightcond1}
E[\Phi(Y^h)]=
{h^2\over c_0}\,E^{(9)}_0\Big[\Phi(R_t,0\leq t\leq h)\,R_h^{-4}\,\exp
\Big(\int_0^h{dt\over R_t^2}\,G\Big({R_t\over
\sqrt{h-t}}\Big)\Big)\Big].
\ee
Suppose that on the same probability space where $Y^h$ is defined, we are
also given two random point measures $\m$ and $\m'$ on $\R_+\times\Omega$,
which conditionally given $Y^h$ are independent Poisson point measures
with intensity
\be
\label{Poisson2}
\mu^h_0(Y^h;dr\,d\omega):=2\;\ind{[0,h]}(r)\;\ind{[0,h-r]}(H(\omega))
\,\ind{\{\r(\omega)\subset]-Y^h_r,\infty[\}}
\,dr\,\N_0(d\omega).
\ee
From the continuity properties of the ``reconstruction mapping'' $\Gamma$,
it should now be clear that the probability measures
$\N^{h,(\varepsilon)}$ converge as $\varepsilon\to 0$ to the measure
$\ov\N_0^h$ defined as the law of $\Gamma(Y^h,\m,\m')$. Here we
leave some easy technical details to the reader.

Let us prove the second assertion of Theorem \ref{heightcond}. Let us 
fix $s_1>0$, and let $\psi$ be a continuous function on $\R_+$ with
compact support contained in $]0,\infty[$. Let $F$ be a bounded
continuous function
on $\Omega$. It follows from Theorem \ref{exc-cond} that
\begin{eqnarray}
\label{heightcond2}
\lim_{\varepsilon\to 0}\varepsilon^{-4}\;
\N_0\Big(\psi(\zeta_{s_1})\,F(W)\,
\ind{\{\un W>-\varepsilon\}}\Big)&=&{2\over 21}\;
\N_0(\sigma^{-1}\psi(\zeta^{[s_*]}_{s_1})\,F(W^{[s_*]}))\nonumber\\
&=&{2\over 21}\,\ov\N_0(\psi(\zeta_{s_1})\,F(W)).
\end{eqnarray}
To see this, apply Theorem \ref{exc-cond} with a function $\varphi$ such
that $s\varphi(s)$ vanishes on a neighborhood of $0$ and is
identically equal to $1$ on $[s_1,\infty[$. 

On the other hand, we have also
\begin{eqnarray}
\label{heightcond3}
\varepsilon^{-4}\;
\N_0\Big(\psi(\zeta_{s_1})F(W)
\ind{\{\un W>-\varepsilon\}}\Big)\!\!&=&
\varepsilon^{-4}\int_0^\infty {dh\over 2h^2}\,
\N^h_0\Big(\psi(\zeta_{s_1})\,F(W)\,
\ind{\{\un W>-\varepsilon\}}\Big)\nonumber\\
&=&\!\!\int_0^\infty\!\! {dh\over 2h^2}\,\varepsilon^{-4}\N^h_0(\un
W>-\varepsilon)\times
\N^{h,(\varepsilon)}_0(\psi(\zeta_{s_1})F(W)).
\end{eqnarray}
We pass to the limit $\varepsilon\to 0$ in the right-hand side
of (\ref{heightcond3}), using (\ref{AbW2}) and the first assertion of the
theorem, which gives
$$\lim_{\varepsilon\to 0} 
\N^{h,(\varepsilon)}_0(\psi(\zeta_{s_1})\,F(W))
=\ov \N^h_0(\psi(\zeta_{s_1})\,F(W)).$$
To justify dominated convergence, first note that
\be
\label{Lgtech}
\varepsilon^{-4}\,\N^h_0(\un
W>-\varepsilon)=\varepsilon^{-4}\,\N^1_0(\un
W>-{\varepsilon\over\sqrt{h}})\leq {C\over h^2}.
\ee
Furthermore, by comparing the intensity measures in
(\ref{Poisson1}) and (\ref{Poisson2}), we get that the
distribution of $\sigma$ under $\N^{h,(\varepsilon)}_0$ is stochastically bounded
by the distribution of $\sigma$ under $\N^h_0$. Hence,
$$\N^{h,(\varepsilon)}_0(\psi(\zeta_{s_1})\,F(W))
\leq C'\N^{h,(\varepsilon)}_0(\sigma>s_1)
\leq C'\N^h_0(\sigma>s_1)\leq C_{(s_1)}\,\exp(-{C'_{(s_1)}\over h^2}),$$
where $C_{(s_1)}$ and $C'_{(s_1)}$ are positive constants depending on
$s_1$.

The previous observations allow us to apply the dominated convergence 
theorem to the right-hand side of (\ref{heightcond3}), and to get
$$\lim_{\varepsilon\da
0}\varepsilon^{-4}\;\N_0\Big(\psi(\zeta_{s_1})\,F(W)
\ind{\{\un W>-\varepsilon\}}\Big)
={c_0\over 2}\int_0^\infty {dh\over
h^4}\,\ov\N^h_0(\psi(\zeta_{s_1})\,F(W)).$$
Comparing with (\ref{heightcond2}) now completes the proof. \cq

At this point we have obtained two distinct descriptions of $\ov\N_0$:
\begin{description}
\item{$\bullet$} The law of $\sigma$ under $\ov\N_0$ has density
$(8\pi)^{-1/2}s^{-5/2}$, and the conditional distribution
$\ov\N_0(\cdot\mid \sigma=s)$ is the law under $\N^{(s)}_0$
of the re-rooted snake $W^{[s_*]}$. 
\item{$\bullet$} The law of $H$ under $\ov\N_0$ has density ${21c_0\over
4}
\,h^{-4}$, and the conditional distribution $\ov\N_0(\cdot\mid H=h)$
can be reconstructed from the ``spine'' $Y^h$ and the Poisson
point measures $\m$ and $\m'$ as explained in the
proof of Theorem \ref{heightcond}.
\end{description}
If we think of analogous results for the It\^o measure of Brownian
excursions, it is tempting to look for a more Markovian description of
$\ov\N_0$. It is relatively easy to see that the process 
$((\zeta_s,W_s),s>0)$ is Markovian under $\ov\N_0$, and to describe its
transition kernels (informally, this is the Brownian snake conditioned
not to exit $]0,\infty[$ -- compare with \cite{AS}). One would then like
to have an explicit formula for entrance laws, that is for the law of
of $(\zeta_s,W_s)$ under $\ov\N_0$, for each fixed $s>0$. Such explicit
expressions seem difficult to obtain. See however the calculations
in Section 5.

In the final part of this section, we investigate the limiting 
behavior of the measures $\ov\N_0(\cdot\mid H=h)$ as $h\to \infty$. This
leads to a (one-dimensional) Brownian snake
conditioned to stay positive and to live forever. The motivation 
for introducing such a process comes from the fact that it is expected 
to appear in scaling limits of discrete trees coding random 
quadrangulations: See the recent work of Chassaing and Durhuus \cite{CD}.

Before stating our result, we give a description of the limiting process.
Let $Z=(Z_t,t\geq 0)$ be a nine-dimensional Bessel process started at
$0$. Conditionally given $Z$, let
$$\p=\sum_{i\in I} \delta_{(h_i,\omega_i)}$$
be a Poisson point measure on $\R_+\times \Omega$ with intensity
$$2\;\ind{\{\r(\omega)\subset]-Z_r,\infty[\}}\,dr\,\N_0(d\omega).$$
We may and will assume that $\p$
is constructed in the following way. Start from a Poisson point measure
$$\q=\sum_{j\in J}\delta_{(\ov h_j,\ov\omega_j)}$$
with intensity
$2\;dr\,\N_0(d\omega)$, and assume that $\q$ is independent of $Z$. Then
set
$$\p=\sum_{j\in J} \ind{\{\r(\ov \omega_j)\subset ]-Z_{\ov h_j},\infty[\}}
\,\delta_{(\ov h_j,\ov\omega_j)}.$$

We then construct our conditioned snake $W^\infty$
from the pair $(Z,\p)$. This is very similar to the reconstruction
mapping that was already used in the proof of Theorem \ref{heightcond}.
To simplify notation, we put
$$\sigma_i=\sigma(\omega_i)\;,\quad\zeta^i_s=\zeta_s(\omega_i)\;,\quad
W^i_s=W_s(\omega_i)$$
for every $i\in I$ and $s\geq 0$. For every $\ell\in[0,h]$, we set
$$\tau_\ell=\sum_{i\in I} \ind{\{h_i\leq \ell\}}\,\sigma_i.$$
Then, if $s\geq 0$, there is a unique $\ell$ such that
$\tau_{\ell-}\leq s\leq \tau_\ell$, and:
\begin{description}
\item{$\bullet$} Either there is a (unique) $i\in I$ such that $\ell=h_i$,
and we set
\ba
&&\zeta^\infty_s=\ell+\zeta^i_{s-\tau_{\ell-}}\;,\\
&&W^\infty_s(t)=\left\{\begin{array}{ll} Z_t&{\rm if}\ t\leq
\ell\;,\\ Z_\ell+W^i_{s-\tau_{\ell-}}(t-\ell)\qquad&{\rm if}\
\ell<t\leq
\zeta^\infty_s\;;
\end{array}
\right.
\ea
\item{$\bullet$} Or there is no such $i$, and we set
\ba
&&\zeta^\infty_s=\ell\;,\\
&&W^\infty_s(t)=Z_t\;,\qquad t\leq \ell\;.
\ea
\end{description}
It is easy to verify that these prescriptions define a 
continuous process $W^\infty$ with values in
$\W$. We denote by $\ov\N^\infty_0$ the law of $W^\infty$.

\begin{theorem}
\label{nonextinct}
The probability measures $\ov\N^h_0$ converge 
to $\ov\N^\infty_0$ when $h\to\infty$. 
\end{theorem}
\proof We rely on the explicit description of
$\ov\N^h_0$ obtained in the proof of Theorem \ref{heightcond}. Let 
$Y^h=(Y^h_t,0\leq t\leq h)$ be as in (\ref{heightcond1}).

\begin{lemma}
\label{nonextinctlem}
The processes $(Y^h_{t\wedge h},t\geq 0)$ converge in distribution to $Z$
as
$h\to\infty$.
\end{lemma}
\proof Let $A>0$ and let $\Phi$ be a bounded continuous function
on $C([0,A],\R_+)$. By (\ref{heightcond1}), if $h\geq A$,
$$E[\Phi(Y^h_t,0\leq t\leq A)]=
{h^2\over c_0}\,E^{(9)}_0\Big[\Phi(R_t,0\leq t\leq A)\,R_h^{-4}\,\exp
\Big(\int_0^h{dt\over R_t^2}\,G\Big({R_t\over
\sqrt{h-t}}\Big)\Big)\Big].$$
We apply the Markov property at time $A$ in the
right-hand side, and write $h=A+a$ to simplify notation:
\begin{eqnarray}
\label{nonext1}
&&E[\Phi(Y^h_t,0\leq t\leq A)]=
{h^2\over c_0}\,E^{(9)}_0\Big[\Phi(R_t,0\leq t\leq A)\,
\exp
\Big(\int_0^A{dt\over R_t^2}\,G\Big({R_t\over
\sqrt{h-t}}\Big)\Big)\nonumber\\
&&\hspace{6cm} \times\; E^{(9)}_{R_A}\Big[R_a^{-4}\,\exp
\Big(\int_0^a{dt\over R_t^2}\,G\Big({R_t\over
\sqrt{a-t}}\Big)\Big)\Big]\Big].
\end{eqnarray}
From the bound $0\leq G(x)\leq 2x^2$, it is immediate that
\be
\label{nonext2}
1\leq \exp\Big(\int_0^A {dt\over R_t^2}\,G\Big({R_t\over
\sqrt{h-t}}\Big)\Big)\leq \exp({2A\over a}).
\ee
On the other hand, a scaling argument gives 
$${h^2\over c_0}\,E^{(9)}_{R_A}\Big[R_a^{-4}\,\exp
\Big(\int_0^a{dt\over R_t^2}\,G\Big({R_t\over
\sqrt{a-t}}\Big)\Big)\Big]
=({h\over a})^2c_0^{-1}\,E^{(9)}_{R_A/\sqrt{a}}\Big[R_1^{-4}\,\exp
\Big(\int_0^1{dt\over R_t^2}\,G\Big({R_t\over
\sqrt{1-t}}\Big)\Big)\Big].$$
From (\ref{heighttech3}), we know that
\be
\label{nonext3}
\lim_{x\da 0} E^{(9)}_{x}\Big[R_1^{-4}\,\exp
\Big(\int_0^1{dt\over R_t^2}\,G\Big({R_t\over
\sqrt{1-t}}\Big)\Big)\Big]=E^{(9)}_{0}\Big[R_1^{-4}\,\exp
\Big(\int_0^1{dt\over R_t^2}\,G\Big({R_t\over
\sqrt{1-t}}\Big)\Big)\Big]=c_0.
\ee
We can use (\ref{nonext2}) and (\ref{nonext3}) to pass to the limit $h\to
\infty$
in the right-hand side of (\ref{nonext1}). The justification of
dominated convergence is easy thanks to the bounds we obtained 
when proving (\ref{heighttech3}). It follows that
$$\lim_{h\to\infty}E[\Phi(Y^h_t,0\leq t\leq A)]=
E^{(9)}_0[\Phi(R_t,0\leq t\leq A)]$$
which was the desired result. \cq

We can now complete the proof of Theorem \ref{nonextinct}.
By Lemma \ref{nonextinctlem} and the Skorokhod representation
theorem, we may
assume that
$(Y^h_{t\wedge h})_{t\geq 0}$ converges to $(Z_t)_{t\geq 0}$
uniformly on every compact subset of $\R_+$, a.s.

Recall the description of $\ov\N_0^h$ as the law of $\Gamma(Y^h,\m,\m')$
in
the proof of Theorem \ref{heightcond}: According to this description, we
can construct a process 
$(W^h_s)_{s\leq \alpha_h}$ having the distribution
of $(W_s)_{s\leq \alpha}$ under $\ov\N^h_0$, by the same formulas
we used to define $W^\infty$ from the pair $(Z,\p)$, 
provided that
$Z$ is replaced  by $Y^h$, the point measure $\p$ is replaced by 
$$\n^h:=\sum_{j\in J} 
\ind{\{\r(\ov\omega_j)\subset]-Y^h_{\ov h_j},\infty[\}}\,
\ind{\{\ov h_j<h,\,H(\ov \omega_j)<h-\ov h_j\}}\,
\delta_{(\ov h_j,\ov\omega_j)}$$
(note that the conditional distribution of $\n^h$ knowing $Y^h$ is
that of a Poisson point measure with intensity $\mu^h_0(Y^h;dr\,d\omega)$,
as required) and we restrict our attention to
$$s\leq \alpha_h:=\int \n^h(dr\,d\omega)\,\sigma(\omega).$$
When $h\to \infty$, the constraints $\{\ov h_j<h,\,H(\ov \omega_j)<h-\ov
h_j\}$ in the definition of $\n^h$ play no role, and the convergence
of $Y^h$ to $Z$ implies that $\n^h$ converges to $\p$, in a sense that
can easily be made precise. It is then a straightforward exercise to
verify that
$$\lim_{h\to\infty} \ (W^h_{s\wedge
\alpha_h})_{s\geq 0}=(W^\infty_s)_{s\geq 0}$$
uniformly on every compact subset of $\R_+$, a.s. The 
statement of Theorem \ref{nonextinct} follows. \cq

\section{Finite-dimensional marginal distributions under $\ov\N_0$}

Our goal in this section is to get an analogue of formula (\ref{margi})
when $\N_x$ is replaced by the conditional measure $\ov\N_0$. This result
will be formally analogous to (\ref{margi}) but the role of Brownian motion
for the spatial displacements will be played by the nine-dimensional
Bessel process. More precisely, recall the notation before (\ref{margi}),
and let $x\geq 0$.
For a fixed marked tree $\theta=(\t,(h_u)_{u\in\t})$, $(\xi^u,u\in\t)$
are independent linear Brownian motions
under the probability mesasure $Q^\theta_x$. 
Under the same probability measure, we construct inductively a 
collection of nine-dimensional Bessel processes $(\ov\xi^u,u\in\t)$
by first requiring that $\ov\xi^\varnothing$
is obtained as the solution of the stochastic differential equation
$$\left\{\begin{array}{ll}
d\ov\xi^\varnothing_t=d\xi^\varnothing_t+
{\displaystyle{4\over \ov\xi^\varnothing_t}}\,dt\;,
\qquad&0\leq t\leq h_\varnothing\;,\\
\ov\xi^\varnothing_0=x\;,&
\end{array}
\right.
$$
and then, for every $u\in\t\backslash\{\varnothing\}$, constructing 
$\ov\xi^u$ as the solution of
$$\left\{\begin{array}{ll}
d\ov\xi^u_t=d\xi^u_t+{\displaystyle{4\over \ov\xi^u_t}}\,dt\;,
\qquad&0\leq t\leq h_u\;,\\
\ov\xi^u_0=\ov\xi^{\pi(u)}_{h_{\pi(u)}}.&
\end{array}
\right.
$$
We then define $(\ov V_a, a\in \wt \theta)$ by the formula
$\ov V_{p_\theta(u,\ell)}=\ov\xi^u_\ell$ for every $u\in \t$ and
$\ell\in[0,h_u]$. Finally, for every leaf $a$ of $\wt\theta$, we define the
stopped path $\ov\w^{(a)}$ from $(\ov V_a, a\in \wt \theta)$ in the same
way as $\w^{(a)}$ was defined from $(V_a,a\in\wt \theta)$. Recall the notation
$L(\theta)$ for the set of leaves of $\wt \theta$, and $I(\theta)$ for the set of its nodes.

\begin{theorem} 
\label{marginalscond}
Let $p\geq 1$ be an integer.
Let $F$ be a symmetric nonnegative measurable function on $\W^p$. Then,
$$\ov\N_0\Big(\int_{]0,\sigma[^p} \!\!ds_1\ldots
ds_pF(W_{s_1},\ldots,W_{s_p})\Big) 
=2^{p-1}p!\!\! \int\!\!
\Lambda_p(d\theta)Q^\theta_0\Big[F((\ov\w^{(a)})_{a\in
L(\theta)})\!\!\prod_{a\in I(\theta)}(\ov V_a)^4\!\!
\prod_{a\in L(\theta)} (\ov V_a)^{-4}\Big].
$$
\end{theorem} 

\proof We may assume that $F$ is continuous and bounded above by $1$, and that there
exist positive constants $\delta$ and $M$ such that
$F(\w_1,\ldots,\w_p)=0$ as soon as $\zeta_{(\w_i)}\notin[\delta,M]$
for some $i$. The proof
will be divided in several steps.

{\em Step 1}. To simplify notation, we write $\r(V)=\{V_a,a\in \wt \theta\}$,
for the range of $V$, or equivalently for the union of the ranges of 
$\w^{(a)}$ for $a\in L(\theta)$. We first
apply Theorem
\ref{marginals} to compute 
\ba
&&\hspace{-5mm}\N_0\Big(\int_{]0,\sigma[^p} ds_1\ldots
ds_p\,F(W_{s_1},\ldots,W_{s_p})\;\ind{\{\r\subset]-\varepsilon,\infty[\}}\Big)\\
&&\hspace{-5mm}=
p!2^{p-1}\int\Lambda_p(d\theta)Q_0^\theta\cg F((\w^{(a)})_{a\in
L(\theta)})\ind{\{\r(V)\subset]-\varepsilon,\infty[\}}
\exp\pg-4\int \l_\theta(da)\,
\N_{0}(\r\subset]-\vep-V_a,\infty[)\pd\cd\\ 
&&\hspace{-5mm}=
p!2^{p-1}\int\Lambda_p(d\theta)Q_0^\theta\cg F((\w^{(a)})_{a\in
L(\theta)})\ind{\{\r(V)\subset]-\varepsilon,\infty[\}}
\exp\pg-6\int {\l_\theta(da)\over (V_a+\varepsilon)^2}\pd\cd\\
&&\hspace{-5mm}=
p!2^{p-1}\int\Lambda_p(d\theta)Q_\varepsilon^\theta\cg
F((-\varepsilon+\w^{(a)})_{a\in L(\theta)})\ind{\{\r(V)\subset]0,\infty[\}}
\exp\pg-6\int {\l_\theta(da)\over (V_a)^2}\pd\cd.\\
\ea
We then use Proposition \ref{Bessel} inductively to see that 
\ba
&&\varepsilon^{-4}Q_\varepsilon^\theta\cg
F((-\varepsilon+\w^{(a)})_{a\in L(\theta)})\ind{\{\r(V)\subset]0,\infty[\}}
\exp\pg-6\int {\l_\theta(da)\over (V_a)^2}\pd\cd\\
&&\qquad=Q^\theta_\varepsilon\cg
F((-\varepsilon+\ov\w^{(a)})_{a\in L(\theta)})
\prod_{a\in I(\theta)} (\ov V_a)^4\prod_{a\in L(\theta)}(\ov V_a)^{-4}\pd\cd.
\ea
We have thus proved that
\begin{eqnarray}
\label{Marg1}
&&\vep^{-4}\N_0\Big(\int_{]0,\sigma[^p} ds_1\ldots
ds_p\,F(W_{s_1},\ldots,W_{s_p})\;\ind{\{\un W>-\varepsilon\}}\Big)
\nonumber
\\&&\qquad=p!2^{p-1}\int\Lambda_p(d\theta)
Q^\theta_\varepsilon\cg
F((-\varepsilon+\ov\w^{(a)})_{a\in L(\theta)})
\prod_{a\in I(\theta)} (\ov V_a)^4\prod_{a\in L(\theta)}(\ov V_a)^{-4}\pd\cd.
\end{eqnarray}

{\em Step 2}. We focus on the right-hand side of (\ref{Marg1}). 
Our goal is to prove that 
\begin{eqnarray}
\label{Limitedroite}
&&\lim_{\varepsilon\to0}\int\Lambda_p(d\theta)Q^\theta_\varepsilon\cg
F((-\varepsilon+\ov\w^{(a)})_{a\in L(\theta)})
\prod_{a\in I(\theta)} (\ov V_a)^4\prod_{a\in L(\theta)}(\ov V_a)^{-4}\pd\cd
\nonumber\\
&&\qquad=\int\Lambda_p(d\theta)Q^\theta_0\cg
F((\ov\w^{(a)})_{a\in L(\theta)})
\prod_{a\in I(\theta)} (\ov V_a)^4\prod_{a\in L(\theta)}(\ov V_a)^{-4}\pd\cd.
\end{eqnarray}
We first state a lemma.
\begin{lemma}{\label{Marg2}}
We
have
$$Q^\theta_\varepsilon\cg
\prod_{a\in I(\theta)} (\ov V_a)^4\prod_{a\in L(\theta)}(\ov V_a)^{-4}\pd\cd
\leq E^{(9)}_\varepsilon[R_{D(\theta)}^{-4}]$$ 
where
$D(\theta)=\max\{d_\theta(0,a):a\in L(\theta)\}$.
\end{lemma}

\proof 
We argue by induction on $p$. If $p=1$, the result is immediate, with an equality.
Let $p\geq 2$ and let us assume that the result holds at order 
$1,2,\ldots,p-1$. Let $\theta=(\t,(h_u,u\in\t))$ be a
marked tree with $p$ leaves. Write $h=h_\varnothing$. By decomposing $\theta$ 
at its first branching point, we get two marked trees $\theta'\in\T_j$, 
and $\theta''\in\T_{p-j}$, for some $j\in\{1,\ldots,p-1\}$, in such a way that
\ba
&&Q^\theta_\varepsilon\cg
\prod_{a\in I(\theta)} (\ov V_a)^4\prod_{a\in L(\theta)}(\ov V_a)^{-4}\cd\\
&&\qquad= E_\vep^{(9)}\cg R_h^4 \,Q^{\theta'}_{R_h}\cg
\prod_{a\in I(\theta')} (\ov V_a)^4\prod_{a\in L(\theta')}(\ov V_a)^{-4}\cd
Q^{\theta''}_{R_h}\cg
\prod_{a\in I(\theta'')} (\ov V_a)^4\prod_{a\in L(\theta'')}(\ov V_a)^{-4}\cd\cd\\
&&\qquad\leq
E_\vep^{(9)}\cg R_h^4\;E_{R_h}^{(9)}[R_{D(\theta')}^{-4}]\;
E_{R_h}^{(9)}[R_{D(\theta'')}^{-4}]\cd.
\ea
We have used the induction hypothesis in the last inequality. 
We now observe that $D(\theta)=h+\max\{D(\theta'),D(\theta'')\}$. Assume for
definiteness that $D(\theta)=h+D(\theta')$. Using the bound (\ref{easybound}) and the 
Markov property we get
$$
E_\vep^{(9)}\cg R_h^4\;E_{R_h}^{(9)}[R_{D(\theta')}^{-4}]\;
E_{R_h}^{(9)}[R_{D(\theta'')}^{-4}]\cd
\leq E^{(9)}_\varepsilon\cg E^9_{R_h}
[R_{D(\theta')}^{-4}]\cd=E_\vep^{(9)}[R_{h+D(\theta')}^{-4}]
=E_\vep^{(9)}[R_{D(\theta)}^{-4}].
$$
This completes the proof of the lemma.\cq

As a consequence of Lemma \ref{Marg2}, we get the bound
\ba
Q^\theta_\varepsilon\cg F((-\varepsilon+\ov\w^{(a)})_{a\in L(\theta)})
\prod_{a\in I(\theta)} (\ov V_a)^4\prod_{a\in L(\theta)}(\ov V_a)^{-4}\pd\cd
&\leq& E^{(9)}_\varepsilon[R_{D(\theta)}^{-4}]\prod_{a\in L(\theta)}\ind{\{\delta\leq
d_\theta(0,a)\leq M\}}\\
&\leq&E^{(9)}_0[R_{D(\theta)}^{-4}]\prod_{a\in L(\theta)}\ind{\{\delta\leq
d_\theta(0,a)\leq M\}}\\
&=&{E^{(9)}_0[R_1^{-4}]\over D(\theta)^2}\prod_{a\in L(\theta)}\ind{\{\delta\leq
d_\theta(0,a)\leq M\}}.
\ea
The last quantity is clearly integrable
with respect to the measure $\Lambda_p(d\theta)$. In addition, using the
continuity of $F$, it is easy to verify that
\ba
&&\lim_{\vep\to 0} Q^\theta_\varepsilon\cg
F((-\varepsilon+\ov\w^{(a)})_{a\in L(\theta)})
\prod_{a\in I(\theta)} (\ov V_a)^4\prod_{a\in L(\theta)}(\ov V_a)^{-4}\pd\cd\\
&&\qquad=Q^\theta_0\cg
F((\ov\w^{(a)})_{a\in L(\theta)})
\prod_{a\in I(\theta)} (\ov V_a)^4\prod_{a\in L(\theta)}(\ov V_a)^{-4}\pd\cd
\ea
An application of the dominated convergence theorem now leads to
(\ref{Limitedroite}).

{\em Step 3}. We now consider the left-hand side of formula (\ref{Marg1}). For every
$0<b<1$, we consider the continuous function $\phi_b:\R_+\la[0,1]$ such that $\phi_b(s)=1$ 
for every $s\in[b,1/b]$, $\phi_b(s)=0$ for every $s\in\R_+\backslash]b/2,2/b[$,
and $\phi_b$ is linear on $[b/2,b]$ and on $[1/b,2/b]$. From Theorem
\ref{exc-cond} and the definition of $\ov\N_0$, we get
\begin{eqnarray}
\label{Limitegauche2}
&&\lim_{\vep\rightarrow0}\vep^{-4}\N_0\pg\phi_b(\sg)
\int_{]0,\sigma[^p} ds_1\ldots
ds_p\,F(W_{s_1},\ldots,W_{s_p})\;\ind{\{\un{W}>-\vep\}}\pd\nonumber\\
&&\qquad=
\ov{\N}_0\pg\phi_b(\sg)\int_{]0,\sigma[^p} ds_1\ldots
ds_p\,F(W_{s_1},\ldots,W_{s_p})\pd.
\end{eqnarray}

\begin{lemma}\label{Limitegauche1}
The following convergence holds:
$$\lim_{b\to 0}\sup_{\vep\in(0,1)}\vep^{-4}
\N_0\pg(1-\phi_b(\sg))\int_{]0,\sigma[^p} ds_1\ldots
ds_p\,F(W_{s_1},\ldots,W_{s_p})\ind{\{\un{W}>-\vep\}}\pd=0.$$
\end{lemma}

\proof We first observe that
\begin{eqnarray}
\label{Lg1a}
&&\vep^{-4}\N_0\pg\ind{\{\sg<b\}}\int_{]0,\sigma[^p} ds_1\ldots
ds_p\,F(W_{s_1},\ldots,W_{s_p})\;\ind{\{\un{W}>-\vep\}}\pd\nonumber \\
&&\ \leq b^p \vep^{-4}
\N_0\pg \sup_{s\in[0,\sg]}\zeta_s>\dl\;,\;\un{W}>-\vep\pd
=b^p\varepsilon^{-4}\int_\delta^\infty {dh\over 2h^2}\,\N_0^h(\un{W}>-\vep)
\leq b^p\frac{C}{6\dl^3}
\end{eqnarray}
where the constant $C$ is such that
$\varepsilon^{-4}\N_0^h(\un{W}>-\vep)\leq Ch^{-2}$, for every $h>0$ and
$0<\varepsilon<1$ (cf (\ref{Lgtech})). On the other hand, the Cauchy-Schwarz
inequality gives
\ba
&&\vep^{-4}\N_0\pg\ind{\{\sg>1/b\}}\int_{]0,\sigma[^p} ds_1\ldots
ds_p\,F(W_{s_1},\ldots,W_{s_p})\;\ind{\{\un{W}>-\vep\}}\pd\\
&&\leq\pg\vep^{-4}\N_0\pg\pg\int_{]0,\sigma[^p} ds_1\ldots
ds_p\,F(W_{s_1},\ldots,W_{s_p})\pd^2\;\ind{\{\un{W}>-\vep\}}\pd\pd^{1/2}
\pg\vep^{-4}\N_0\pg\sg>{1\over b},\un{W}>-\vep\pd\pd^{1/2}.\\
\ea
Note that we may write 
\ba
&&\N_0\pg\pg\int_{]0,\sigma[^p} ds_1\ldots
ds_p\,F(W_{s_1},\ldots,W_{s_p})\pd^2\ind{\{\un{W}>-\vep\}}\pd\\
&&\qquad=\N_0\pg\int_{]0,\sigma[^{2p}} ds_1\ldots
ds_{2p}\,G(W_{s_1},\ldots,W_{s_{2p}})\;\ind{\{\un{W}>-\vep\}}\pd,
\ea
where $G$ is a nonnegative symmetric function on $\W^{2p}$, which is also bounded
by $1$. 
As a consequence of (\ref{Marg1}) and Lemma \ref{Marg2}, we then get
\ba
&&\vep^{-4}\N_0\pg\pg\int_{]0,\sigma[^p} ds_1\ldots
ds_p\,F(W_{s_1},\ldots,W_{s_p})\pd^{2}\ind{\{\un{W}>-\vep\}}\pd\\
&&\quad\leq (2p)!\,
2^{2p-1}\,E_0^{(9)}[R_1^{-4}]\int\Lambda_{2p}(d\theta)D(\theta)^{-2}
\prod_{a\in L(\theta)}\ind{\{\dl\leq d_\theta(\varnothing,a)\leq M\}}=C(p,\dl,M)<\infty.
\ea
From Theorem \ref{main-cond} and a simple scaling argument, we have
$$\vep^{-4}\N_0\pg\sg>{1\over b},\,\un{W}>-\vep\pd
=\vep^{-4}\int_{b^{-1}}^\infty {ds\over\sqrt{2\pi s^3}}\,
\N^{(s)}_0(\un W>-\varepsilon)\leq C'\,b^{1/2}.$$
By combining these estimates, we get
\be
\label{Lg1b}
\vep^{-4}\N_0\pg\ind{\{\sg>1/b\}}\int_{]0,\sigma[^p} ds_1\ldots
ds_p\,F(W_{s_1},\ldots,W_{s_p})\;\ind{\{\un{W}>-\vep\}}\pd
\leq (C'C(p,\delta,M))^{1/2}\,b^{1/4}.
\ee
Lemma
\ref{Limitegauche1} follows from (\ref{Lg1a}) and (\ref{Lg1b}).
\cq

We can now complete the proof of Theorem \ref{marginalscond}. First, by monotone
convergence,
$$\lim_{b\to 0}
\ov{\N}_0\pg\phi_b(\sg)\int_{]0,\sigma[^p} ds_1\ldots
ds_p\,F(W_{s_1},\ldots,W_{s_p})\pd
=\ov{\N}_0\pg\int_{]0,\sigma[^p} ds_1\ldots
ds_p\,F(W_{s_1},\ldots,W_{s_p})\pd.$$
From (\ref{Limitegauche2}) and Lemma \ref{Limitegauche1}, it then follows that
$$\lim_{\vep\rightarrow0}\vep^{-4}\N_0\pg
\int_{]0,\sigma[^p} \!\!ds_1\ldots
ds_p\,F(W_{s_1},\ldots,W_{s_p})\ind{\{\un{W}>-\vep\}}\pd
=\ov{\N}_0\pg\int_{]0,\sigma[^p} \!\!ds_1\ldots
ds_pF(W_{s_1},\ldots,W_{s_p})\pd.$$
Combining this with (\ref{Marg1}) and (\ref{Limitedroite}) 
gives Theorem \ref{marginalscond}. \cq

\bigskip
\noindent{\bf Acknowledgement.} The first author wishes to thank Philippe Chassaing for
a stimulating conversation which motivated the present work.


\begin{thebibliography}{99}

\bibitem{AW}
{\sc Abraham, R., Werner, W.} (1997) Avoiding probabilities for Brownian snakes and super-Brownian motion.
{\it Electron. J. Probab.} {\bf 2}
no. 3, 27 pp.

\bibitem{AS}
{\sc Abraham, R., Serlet, L.} (2002) Representations of the Brownian snake with drift. {\it Stochastics and
Stochastics
Reports} {\bf 73}, 287-308. 

\bibitem{Al1} 
{\sc Aldous, D.} (1991) The continuum random tree I. {\it Ann. Probab.} {\bf 19},
1-28.

\bibitem{Al2} 
{\sc Aldous, D.} (1991) The continuum random tree. II. An overview. Stochastic analysis (Durham, 1990), 23-70,
London Math. Soc. Lecture Note Ser. 167. Cambridge Univ. Press, Cambridge, 1991. 

\bibitem{Al3} 
{\sc Aldous, D.} (1993) The continuum random tree III. {\it Ann. Probab.} {\bf 21},
248-289.

\bibitem{BDG1}
{\sc Bouttier, J., Di Francesco, P., Guitter, E.} (2003) Random trees between two walls: 
exact partition function.
{\it J. Phys. A} {\bf 36}, 12349-12366.

\bibitem{BDG2}
{\sc Bouttier, J., Di Francesco, P., Guitter, E.} (2003) Statistics of planar graphs viewed from a vertex: a study
via labeled trees. {\it Nuclear Phys. B} {\bf 675}, 631-660.

\bibitem{CD}
{\sc Chassaing, P., Durhuus, B.} (2003)
Statistical Hausdorff dimension of labelled trees and quadrangulations. Preprint.

\bibitem{CS}
{\sc Chassaing, P., Schaeffer, G.} (2004) 
Random planar lattices and integrated superBrownian excursion.
{\it Probab. Th. Rel. Fields} {\bf 128}, 161-212.

\bibitem{Delmas} {\sc Delmas, J.F.} (2003) 
Computation of moments for the length of the one dimensional ISE support. {\it Electron. J. Probab.} {\bf 8}
no. 17, 15 pp. 

\bibitem{Sl1}
{\sc Derbez, E., Slade, G.} (1998) The scaling limit of lattice trees in high dimensions. 
{\it Comm. Math. Phys.} {\bf 198},
69-104.

\bibitem{Duq} {\sc Duquesne, T.} (2003) A limit theorem for the contour process of conditioned 
Galton-Watson trees. {\it Ann. Probab.} {\bf 31}, 996-1027.
    

\bibitem{DuLG} {\sc Duquesne, T., Le Gall, J.F.} (2004) Probabilistic and fractal aspects of L\'evy trees.
{\it Probab. Th. Rel. Fields}, to appear.


\bibitem{EPW} {\sc Evans, S.N., Pitman, J.W., Winter, A.} (2003) Rayleigh processes, real trees
and root growth with re-grafting. {\it Probab. Th. Rel. Fields}, to appear.

\bibitem{Sl2}
{\sc Hara, T., Slade, G.} (2000) The scaling limit of the incipient infinite cluster in high-dimensional
percolation. II. Integrated super-Brownian excursion. 
Probabilistic techniques in equilibrium and nonequilibrium statistical physics. 
{\it J. Math. Phys.} {\bf 41} (2000), 1244-1293. 

\bibitem{Sl3}
{\sc van der Hofstad, R., Slade, G.} (2003) Convergence of critical oriented percolation
to super-Brownian motion above $4+1$ dimensions. {\it Ann. Inst. H. Poincar\'e Probab. Statist.}
{\bf 20}, 413-485. 


\bibitem{JM} {\sc Janson, S., Marckert, J.F.} (2003)
Convergence of discrete snakes. Preprint. 

\bibitem{JR} {\sc Jansons, K.M., Rogers, L.C.G.} (1992) Decomposing the branching Brownian path. {\it Ann. Probab.}
{\bf 2}, 973-986.

\bibitem{Kal} {\sc Kallenberg, O.} (1975) {\it Random Measures}. Academic Press, London.

\bibitem{LG0}
{\sc Le Gall, J.F.} (1991)
Brownian excursions, trees and measure-valued branching
processes. {\sl Ann. Probab.} {\bf 19}, 1399-1439.

\bibitem{LG1}
{\sc Le Gall, J.F.} (1993)
The uniform random tree in a Brownian excursion.
{\sl Probab. Th. Rel. Fields} {\bf 96}, 369-383.

\bibitem{LG2}
{\sc Le Gall, J.F.} (1995)
The Brownian snake and solutions of $\Delta u=u^2$ in a domain. {\sl Probab. Th. Rel.
Fields} {\bf 102}, 393-432.

\bibitem{Zu} {\sc Le Gall, J.F.} (1999) {\it Spatial Branching Processes, Random Snakes and 
Partial Differential Equations}. {\it Lectures in Mathematics ETH Z\"urich}. Birkh\"auser, Boston.

\bibitem{LG}
{\sc Le Gall, J.F.} (2004)
An invariance principle for conditioned Brownian trees. In preparation.

\bibitem{MM2} {\sc Marckert, J.F., A. Mokkadem} (2004) State spaces of the snake and its tour - Convergence
of the discrete snake. {\it J. Theoret. Probability} {\bf 16}, 1015-1046.

\bibitem{MM} {\sc Marckert, J.F., A. Mokkadem} (2004) Limits of normalized quadrangulations. The Brownian
map. Preprint.

\bibitem{RY} {\sc Revuz, D., Yor, M.} (1991) {\it Continuous Martingales and Brownian Motion}.
Springer, Berlin-Heidelberg-New York.

\bibitem{Verwaat} {\sc Verwaat, W.} (1982) A relation between Brownian bridge and Brownian excursion.
{\it Ann. Probab.} {\bf 10}, 234-239.

\bibitem{Yor}
{\sc Yor, M.} (1980) Loi de l'indice du lacet brownien, et distribution de Hartman-Watson.
{\it Z. Wahrsch. verw. Gebiete} {\bf 53}, 71-95.



\end{thebibliography}
\end{document}